\RequirePackage{fix-cm}
\documentclass[smallextended,hidelinks]{svjour3}       
\smartqed  
\usepackage{graphicx,amssymb}

\usepackage{stmaryrd}
\usepackage{xcolor}
\usepackage{algorithm}
\usepackage{algpseudocode} 
\usepackage{placeins}

\usepackage{lineno,hyperref}
\usepackage{mathtools}

\usepackage{enumitem}
\usepackage{tikz}
\usepackage{tcolorbox}
\usepackage{rotating}
\DeclareMathOperator*{\ddiv}{div}
\DeclareMathOperator*{\curl}{curl}
\DeclareMathOperator*{\Curl}{Curl}
\DeclareMathOperator*{\argmin}{arg\,min}
\renewcommand{\div}{{\textrm{ div }}}

\newcommand{\pw}{\mathrlap{\pd}{w}}
\renewcommand{\pw}{{\rm pw}}
\newcommand{\divpw}{{\ddiv}_{\rm pw}}
\newcommand{\apw}{a_{\pw}}

\newcommand{\R}{\mathbb R}

\newcommand{\vF}{v_{\mathcal F}}
\newcommand{\cF}{{\mathcal F}}
\newcommand{\vT}{v_{\mathcal T}}
\newcommand{\vTF}{(\vT,\vF)}

\newcommand{\bn}{\nu} 

\newcommand{\npw}[1]{\trb{ #1 }_{\pw}}

\newcommand{\nablapw}{\nabla_{\pw}}

\newcommand{\Hdivz}{H_0(\ddiv,\oz)}

\newcommand{\bN}{\mathbb N}

\newcommand{\cT}{\mathcal T}
\newcommand{\RT}{\rm RT}
\newcommand{\trb}[1]{|\!|\!|#1|\!|\!|}

\newcommand{\part}{\mathcal T}

\newcommand{\facets}{\mathcal F}
\newcommand{\facetsI}{\mathcal F(\Omega)}

\newcommand{\dPT}[1]{ P_{#1}(T)}

\newcommand{\Vh}{ V_h}

\newcommand{\cR}{ R}

\newcommand{\N}{{\mathbb N}}

\newcommand{\F}{\mathcal{F}}

\newcommand{\T}{\mathcal{T}}

\newcommand{\ENorm}[1]{{\left\vert\kern-0.25ex\left\vert\kern-0.25ex\left\vert #1 
    \right\vert\kern-0.25ex\right\vert\kern-0.25ex\right\vert}}
\newcommand{\capx}{c_{\rm apx}}
\newcommand{\Ca}{C_{\rm 1}}
\newcommand{\Cb}{C_{\rm 2}}

\newcommand{\cstb}{C_{\rm st}}
\newcommand{\ndof}{{\rm ndof}}
\newcommand{\jj}{j_{1,1}}
\newcommand{\CeqT}{\Ca}
\newcommand{\CeqF}{C_{2}}
\usepackage{ulem}

\newcommand{\etares}{\eta_{\rm res}}
\newcommand{\etaHHO}{\eta_{\rm HHO}}

\newcommand{\etaeqnew}{\eta_{{\rm eq}, p}}
\newcommand{\etaeqnewa}[1]{\eta_{{\rm eq}, #1}}

\newcommand{\SpTF}{S_{TF}}
\newcommand{\CH}{C_{\rm H}}
\newcommand{\Ctr}{C_{\rm tr}}
\newcommand{\Ceff}{\cnst{cnst:res_eff}}

\newcommand{\osck}{\mathrm{osc}_k}
\newcommand{\osc}[1]{\mathrm{osc}_{#1}}
\newcommand{\qRT}{q_{RT}}
\newcommand{\oF}{\omega(F)}
\newcommand{\QQ}{Q}

\newcommand{\Rzh}{\QQ_{z,h}}
\newcommand{\RzhD}{\Rzh^\Delta}

\newcommand{\GG}{G_h}
\newcommand{\QtoG}{{G}}
\newcommand{\oz}{\omega(z)}

\newcounter{cnt}
\newcommand{\newcnst}{%
	\refstepcounter{cnt}%
	\ensuremath{C_{\thecnt}}}
\newcommand{\cnst}[1]{\ensuremath{C_{\ref{#1}}}}

\begin{document}
\title{Stabilization-free HHO a posteriori error control 
\thanks{
This work has been supported by the Deutsche Forschungsgemeinschaft (DFG) in the Priority Program 1748 
{\it{Reliable simulation techniques in solid mechanics. Development of non-standard discretization methods, mechanical and mathematical analysis}} under the projects BE 6511/1-1  and CA 151/22-2 as well as the European Union’s Horizon 2020 research and innovation
programme (Grant agreement No. 891734). The third author is also supported by the Berlin Mathematical School.
}
}


\author{Fleurianne Bertrand \and Carsten~Carstensen \and Benedikt Gr\"a\ss le \and Ngoc Tien Tran}

\authorrunning{Bertrand, Carstensen, Gr\"a\ss le, Tran} 

\institute{Fleurianne~Bertrand \at
              University of Twente \\
              \email{f.bertrand@utwente.nl}           
           \and
           Carsten~Carstensen 
           \and 
           Benedikt~Gr\"a\ss le  \at
              Humboldt-Universit\"at zu Berlin, Germany\\
              \email{cc, graesslb@math.hu-berlin.de}
            \and
            Tran Ngoc Tien \at Friedrich-Schiller Universit\"at Jena\\
            \email{ngoc.tien.tran@uni-jena.de}
}

\date{}

\maketitle

\begin{abstract}
The known a~posteriori error analysis of hybrid high-order methods (HHO) treats the stabilization contribution as part of the error and as part of the error estimator for an efficient and reliable error control. 
This paper circumvents the stabilization contribution on simplicial meshes  and arrives at a stabilization-free error analysis with an explicit residual-based a~posteriori error estimator for adaptive mesh-refining as well as an equilibrium-based guaranteed upper error bound (GUB). 
Numerical evidence in a Poisson model problem supports that the GUB leads to realistic upper bounds for the displacement error in the piecewise energy norm.
The adaptive  mesh-refining algorithm associated to the explicit residual-based a posteriori error estimator recovers the optimal convergence rates in computational benchmarks.
\keywords{hybrid high-order  \and a posteriori 
\and 
guaranteed upper error bounds
\and adaptive mesh refinement \and equilibration \and stabilization-free
\and
computational comparisons}
\end{abstract}

\section{Introduction}
\label{intro}
Hybrid high-order methods (HHO) were introduced in \cite{Di-Pietro.Ern:15,Di-Pietro.Ern.ea:14} and are examined in the textbooks \cite{DiPietroDroniou2020,ern_finite_2021-2}
as a promising class of flexible nonconforming discretization methods for partial differential equations 
that involve a  parameter-free stabilization term for the link between the volume and skeletal variables. 
\subsection{Known a posteriori error estimator}
The a priori error analysis of HHO involves the stability terms in extended norms as part of the methodology and motivated a first explicit residual-based 
a~posteriori error estimator in \cite{DiPietroDroniou2020} with a reformulation of the stabilization in the upper bound.
Let $s_h(u_h,u_h)$ denote the stabilization at the discrete solution $u_h\in V_h$ and let the (elliptic) reconstruction $Ru_h$ of $u_h$
denote a piecewise polynomial of degree at most $k+1$ that approximates $u\in H^1(\Omega)$, cf.\ \eqref{eqn:s_h0} and Section \ref{sec:residual} below for further details.
Then a possible error term reads 
\begin{align}\label{eqn:tot_error}  
    \textup{total error}^2:=\| \nabla_\textup{pw}(u-Ru_h)\|_{L^2(\Omega)}^2 +s_h(u_h,u_h).
\end{align}
It is disputable if $s_h(u_h,u_h)\geq0$ is an error contribution, but
if the total error includes $s_h(u_h,u_h) $ (or an equivalent form), then the error estimator 
may also include this term (or a computable equivalent) for a reliable and efficient a~posteriori error control. 
Amongst the many skeletal schemes like (nonconforming) virtual elements, hybridized (weak) discontinuous Galerkin schemes et al., the HHO methodology has a clear and efficacious stabilization
\begin{align}\label{eqn:s_h0}
    s_h(v_h, w_h) \coloneqq \sum_{T\in\T}\sum_{F\in\F(T)} h_F^{-1}\langle \SpTF v_h, \SpTF w_h\rangle_{L^2(F)}
\end{align}
with the abbreviation
$
\SpTF v_h 
	\coloneqq \Pi_{F,k}
	\left(
 v_{\T} +
 (1-\Pi_{T,k})R  v_{h} 
 \right)|_{T}-v_{\F}|_{F}
$
for $v_h=(v_\T, v_\F)\in V_h$ in terms of the $L^2$ projections $\Pi_{K,k}$ onto polynomials of degree at most $k$ on a facet or simplex $K\in\F\cup\T$ of diameter $h_K = \mathrm{diam}(K)$; cf. Subsection \ref{sub:notation} for further details.
The original residual-based estimator $\etaHHO$ from the textbook \cite{DiPietroDroniou2020} for the Poisson model problem $-\Delta u = f$ includes \eqref{eqn:s_h0} and an interpolation $\mathcal{A} Ru_h\in V$ of $Ru_h$ by nodal averaging in
\begin{align*}
	\eta_{\rm HHO}^2 =\ & \| h_{\T}(1 - \Pi_0)(f+{\Delta}_{\rm pw} Ru_h)\|_{L^2(\Omega)}^2
	+ \|\nablapw (1-\mathcal{A})Ru_h\|_{L^2(\Omega)}^2\\
	&+ s_h(u_h, u_h).
\end{align*}
(Multiplicative constants are undisplayed in this introduction for simplicity.)
The results from Theorem 4.3 and 4.7 in \cite{DiPietroDroniou2020} show reliability and efficiency for the total error \eqref{eqn:tot_error} and piecewise polynomial source terms $f\in P_{k+1}(\T)$,
$$\text{total error}^2 \approx \etaHHO^2.$$
\subsection{Stabilization-free a posteriori error control}
There are objections against
the double role of $s_h(u_h,u_h)$ on both sides of the efficiency and reliability estimate. First, the term  $s_h(u,u_h) $  may dominate both sides of the error estimate.
In other words, the total error might be equivalent to $s_h(u_h,u_h)$, 
but the quantity of interest may exclusively be
\[
\textup{error}^2:=\| \nabla_\textup{pw}(u-Ru_h)\|_{L^2(\Omega)}^2.
\]
Second, since the stabilization \eqref{eqn:s_h0} incorporates a negative power of the mesh-size, a reduction property for local refinements remains unclear but is inevitable in the proofs of optimal convergence of an adaptive algorithm \cite{bertrand_opt,carstensen_axioms_2014}.
This paper, therefore, asks a different question about the control of the $\textup{error}$
without the stabilization term \eqref{eqn:s_h0} in the upper bound and introduces two stabilization-free error estimators (multiplicative constants are undisplayed)
\begin{align*}
    \eta_{\text{res}}^2=\ & \|h_\T(  f+{\Delta}_{\rm pw} Ru_h)\|_{L^2(\Omega)}^2
    +\sum_{F\in \cF} h_F\|    [ \nablapw Ru_h ]_F \|_{L^2(F)}^2,\\
    \eta_{\text{eq},p}^2=\ & \osc{k+p}^2(f, \T) + \| Q_p - \nabla_\pw R u_h \|_{L^2(\Omega)}^2 + \|\nablapw (1-\mathcal{A})Ru_h\|_{L^2(\Omega)}^2
\end{align*}
for some parameter $p \in \mathbb{N}_0$.
The explicit residual-based a posteriori error estimator $\etares$ follows from the a posteriori methodology in the spirit of \cite{b9ffd40a,c6537ecf,5be62542,normOfdGrad4Hdiv2015ccdpas} with a piecewise volume residual $f+\Delta_\pw Ru_h$ and the jumps $[\nablapw Ru_h]_F$ across a facet $F$ (on the boundary this is only the tangential component of $\nablapw Ru_h$). 
The equilibrated error estimator $\eta_{\mathrm{eq},p}$ includes the post-processed quantity $Q_p \in RT_{k+p}(\T)$ in the space $ RT_{k+p}(\T)$ of Raviart-Thomas functions of degree $k+p$ for $p \in \mathbb{N}_0$ and the nodal average $\mathcal{A} R u_h \in S^{k+1}_0(\mathcal{T})$ of $R u_h$.
The main results establish reliability and efficiency 
$$\text{error}^2 + \osc{k-1}^2(f, \T) \lesssim \etares^2 \approx \eta_{\mathrm{eq},p}^2\lesssim \text{error}^2 + \osc{q}^2(f, \T)$$
for any $p, q\in \mathbb{N}_0$
up to data-oscillations $\osc{q}^2(f, \T)\coloneqq \|h_\T (1-\Pi_{q})f\|^2_{L^2(\Omega)}$ of the source term $f\in L^2(\Omega)$ and {\em without any stabilization terms}.
Computational benchmarks with adaptive mesh-refinement driven by any of these estimators provide numerical evidence for optimal convergence rates.

\subsection{Further contributions and outline}
The higher-order Crouzeix-Raviart finite element schemes are complicated at least in 3D \cite{Ciarlet2018} and then the HHO methodology is an attractive 
alternative even for simplicial triangulations with partly unexpected advantages like the computation of higher-order guaranteed eigenvalue bounds 
\cite{CEP21}. Higher convergence rates rely on an appropriate adaptive mesh-refining algorithm and hence stabilization-free 
a posteriori error estimators are of particular interest. The recent paper \cite{daveiga2021adaptive} establishes the latter for virtual elements with an over-penalization strategy 
as an extension of \cite{bonito_quasi-optimal_2010} for the discontinuous Galerkin schemes. A  disadvantage is the quantification of  the restriction on the stabilization parameter in
practise and poor  condition for larger parameters. The stabilization-free a posteriori error control in this paper is based on two observations for the HHO schemes on simplicial triangulations.
First, the $P_1$-conforming finite element functions let the stabilization vanish and, second, the divergence-free lowest-order Raviart-Thomas functions are $L^2$ perpendicular to the piecewise gradients $\nablapw R u_h$. In fact, those two fairly general properties lead in Section 2 to a reliable explicit residual-based  a posteriori error estimator. 
In contrast to the simplified introduction above, the paper also focuses on multiplicative constants that lead to the GUB
$$\textup{error}\leq \etares \quad\text{and}\quad \textup{error}\leq \eta_{\mathrm{eq},p};$$
cf. Table \ref{tab:expl_const} for explicit quantities and Theorem \ref{thm:reliability} and Theorem \ref{thm:GUB-HHO} for further details.
\begin{table}
    \centering
    \begin{tabular}{c|ccc}
        $\omega_{\rm max}$ &  $\pi$ & $3\pi/2$ & $2\pi$\\\hline
        $M_{\rm bd}$ & 4 & 6 & 8\\
        $\capx$ & 2.9568 & 6.4642 & 11.3771\\\hline
        $\cstb$ & 26.0893 & 55.8498 & 97.5374\\
        $\newcnst\label{cnst:C1}$ & 2.9718 & 6.4710 & 11.3810\\
        $\newcnst\label{cnst:C2}$ & 7.0495 & 15.2341 & 26.7317
    \end{tabular}
    \caption{Explicit constants $\cnst{cnst:C1}, \cstb, \cnst{cnst:C2}$ for right-isosceles triangles with respect to the maximum interior angle $\omega_{\rm max}$ of the polygonal domain $\Omega$.}
    \label{tab:expl_const}
\end{table}
Numerical comparisons of $\etaHHO$ with $\etares$ and $\eta_{\mathrm{eq},p}$ favour the latter.
Section \ref{sec:blocks} identifies general building blocks of the a posteriori error analysis for discontinuous schemes with emphasis on explicit constants.
An application to HHO leads to the new stabilization-free residual-based estimator $\etares$
in Section \ref{sec:residual}.
The alternative stabilization-free error estimator $\eta_{\mathrm{eq},p}$ follows from an equilibration strategy plus post-processing in Section \ref{sec:eq}.
This paper also contributes to the HHO literature a local equivalence of two stabilizations and the efficiency of the stabilization terms up to data-oscillations in extension of \cite{ErnZanotti2020}.
Numerical comparisons of the different error estimators and an error estimator competition for guaranteed error control of the piecewise energy norm in 2D conclude this paper in Section \ref{sec:Numerical results}.
Three computational benchmarks provide striking numerical evidences for the optimality of the associated adaptive algorithms.
The appendix provides algorithmic details on the computation of the post-processed contribution $\|Q_p-\nablapw Ru_h\|_{L^2(\Omega)}$ in $\etaeqnew$.

\subsection{Overall notation}\label{sub:notation}
Standard notation for Sobolev and Lebesgue spaces and norms apply with $\|\bullet\|\coloneqq\|\bullet\|_{L^2(\Omega)}$ and $\trb{\bullet}\coloneqq\|\nabla\bullet\|_{L^2(\Omega)}$.
In particular,
$H(\ddiv, \Omega)$ is the space of Sobolev functions with weak divergence in $L^2(\Omega)$ and $H(\ddiv=0, \Omega)$ contains only divergence-free functions in $H(\ddiv, \Omega)$.
Throughout this paper, $\T$ denotes a shape-regular triangulation  of the polyhedral bounded Lipschitz domain $\Omega\subset \R^n$ into $n$-simplices with facets $\F$ (edges for $n=2$ and faces for $n=3$) and vertices $\mathcal V$.
Let $\F(\Omega)$ (resp.~$\mathcal{V}(\Omega)$) denote the set of interior facets (resp.~vertices) and $\mathcal{F}(\partial \Omega) \coloneqq \mathcal{F} \setminus \mathcal{F}(\partial \Omega)$ (resp.~$V(\partial \Omega) \coloneqq \mathcal{V} \setminus \mathcal{V}(\Omega)$).
Given $v \in \Omega \to \mathbb{R}^n$ and $w : \Omega \to \mathbb{R}^{2n-3}$, let $\curl v \coloneqq \partial_1 v_2 - \partial_2 v_1$ and $\Curl w \coloneqq (\partial_2 w, -\partial_1 w)^t$ if $n = 2$ and $\curl v \coloneqq (\partial_2 v_3 - \partial_3 v_2, \partial_3 v_1 - \partial_1 v_3, \partial_1 v_2 - \partial_2 v_1)^t$ and $\Curl w \coloneqq \curl w$ if $n = 3$.
For $s\in\R$, let $H^s(\T)$, $H(\ddiv,\T)$, and $H(\curl,\T)$ denote the space of piecewise Sobolev functions with restriction to $T \in \T$ in $H^s(T)$, $H(\ddiv,T)$, and $H(\curl,T)$.
To simplify notation, $H^s(K)$ abbreviates $H^s(\mathrm{int}(K))$ for the open interior $\mathrm{int}(K)$ of a compact set $K$.
The $L^2$-scalar product reads \((\bullet, \bullet)_{L^{2}(\omega)}\) for volumes $\omega\subseteq\Omega$ and \(\langle\bullet, \bullet\rangle_{L^2(\gamma)}\) for surfaces $\gamma\subset \overline{\Omega}$ of co-dimension one; the same symbol applies to scalars and to vectors.
For $V \coloneqq H^1_0(\Omega)$ and $V^* \coloneqq H^{-1}(\Omega)$, let $\langle\bullet, \bullet\rangle$ denote the duality-brackets in $V^* \times V$ for the dual space $V^*$ of $V$ equipped with the operator norm 
$\trb{F}_{*}\coloneqq \sup_{v\in V\setminus\{0\}}|F v|/\trb{v}$ for $F\in V^*$.

Define the energy scalar product $a(v,w) \coloneqq (\nabla v, \nabla w)_{L^2(\Omega)}$ for $v, w\in H^1(\Omega)$ and its piecewise version $\apw(v, w)=(\nablapw v,\nablapw w)_{L^2(\Omega)}$ for $v, w\in H^1(\T)$.
The latter induces the seminorm $\trb{\bullet}_\pw \coloneqq a_\pw(\bullet,\bullet)^{1/2}$ in $H^1(\mathcal{T})$.
Here and throughout the paper, $\nablapw$, $\divpw$, $\curl_\pw$, ${\Delta}_{\rm pw}$, denote the piecewise evaluation of the differential operators $\nabla$, $\ddiv$, $\curl_\pw$, $\Delta$ without explicit reference to the underlying shape-regular triangulation $\T$.

The vector space $P_k(K)$ of polynomials of degrees at most $k\in\N_0$ over a facet or simplex $K\in\F\cup\T$ defines the piecewise polynomial spaces
\begin{align*}
P_k(\T) &\coloneqq\{ p \in L^2(\Omega) \ :\ p_{|T} \in \dPT{k} \text{ for all } T \in \T \}, \\
	P_k(\F) &\coloneqq
		\{ p \in L^2(\facets) \ :\ p_{|F} \in P_k(F) \text{ for all } F \in \F
		\} \ 
\end{align*}
and the space of piecewise Raviart-Thomas functions
\begin{align*}
    RT_{k}^\pw(\T) \coloneqq P_k(\T;\R^n) + xP_k(\T).
\end{align*}
The associated $L^2$ projections read 
$\Pi_{K,k}:L^2(\Omega)\to P_k(K), \Pi_k : L^2(\Omega) \rightarrow P_k(\T),$ and $\Pi_{\facets, k} : L^2(\Omega) \rightarrow P_k(\F)$ with the convention $\Pi_{-1}\coloneqq0$.
Abbreviate $S^{k+1}_0(\T) \coloneqq P_{k+1}(\T) \cap V$ and $RT_k(\T) \coloneqq RT^\pw_k(\T) \cap H(\ddiv,\Omega)$ for all $k \in \mathbb{N}_0$.
The piecewise constant mesh-size function $h_\T\in P_0(\T)$ satisfies $h_{\T|T}\coloneqq h_T$ for $T\in\T$ with the diameter $h_K \coloneqq \mathrm{diam}(K)\in P_0(K)$ of $K\in \F\cup\T$.

If not explicitly stated otherwise, constants are independent of the mesh-size in the triangulation but may depend on the shape-regularity and on the polynomial degree $k$.
The abbreviation  $A\lesssim B$ hides a generic constant $C$  (independent of the mesh-size) in $A\leq C\; B$; $A\approx B$ abbreviates $A\lesssim B\lesssim A$.

\section{Foundations of the a~posteriori error analysis}
\label{sec:blocks}
This section investigates general building blocks of the a posteriori error analysis and revisits arguments from \cite{b9ffd40a,c6537ecf,5be62542,normOfdGrad4Hdiv2015ccdpas} with emphasis on multiply connected domains $\Omega\subset \R^n$ for $n=2,3$.
The general setting of this section results in  reliability for an error estimator that is applicable beyond the HHO methodology.
Consider the weak solution $u\in V = H^1_0(\Omega)$ to the Poisson model problem $-\Delta u=f$ a.e.~in $\Omega$ and $u = 0$ on $\partial\Omega$ for a given source $f\in L^2(\Omega)$; i.e., $u\in V$ satisfies 
\begin{align}
\label{eq:poisson}
a(u,v)= (f,v)_{L^2(\Omega)} \quad \text{for all } v\in V .
\end{align}
An approximation $\QtoG \in L^2(\Omega;\R^n)$ of the gradient $\nabla u\in H(\ddiv, \Omega)$  gives rise to the residual $f+\ddiv \QtoG\in V^* =  H^{-1}(\Omega)$ seen as a linear functional  on $V$, i.e.,
\[
\langle f+\ddiv \QtoG,\varphi\rangle := (f,\varphi)_{L^2(\Omega)} - (\QtoG, \nabla \varphi)_{L^2(\Omega)}
\quad\text{for all }\varphi\in V.
\]
Let $\nu_T$ denote the unit outer normal along the boundary $\partial T$ of each simplex $T\in\T$ and fix the orientation of a unit normal $\nu_F=\pm \nu_T$ for each facet $F\in\F(T)$ of $T$ such
that it matches the outer unit normal $\nu$ of $\partial \Omega$ at the boundary.
The jump $[\QtoG]_F$ of a piecewise function in $m\in\N$ components $\QtoG\in H^1(\T;\R^m)$ reads $[\QtoG]_F\coloneqq \QtoG_{|T_+}-\QtoG_{|T_-}$ on interior facets $F=T_+\cap T_-\in\F(\Omega)$ (with $T_\pm$ labelled such that $\nu_{T_+|F} = \nu_F = -\nu_{T_-|F})$ and $[\QtoG]_F\coloneqq \QtoG$ on the boundary $F\in\F(\partial \Omega)$.
The main result of this section establishes the residual-based error estimator
\begin{align}
\eta^2(\T, \QtoG) \coloneqq & \left(\cnst{cnst:C1}\|h_\T(  f+{\ddiv}_{\rm pw} \QtoG)\|
+\cnst{cnst:C2}\sqrt{\sum_{F\in \cF(\Omega)} \ell(F)\|    [ \QtoG ]_F\cdot \nu_F \|_{L^2(F)}^2}\right)^2\label{eq:eta}\\
&+\CH^2\left(\cnst{cnst:C1}\| h_\T{\curl}_{\rm pw} \QtoG \|
+\cnst{cnst:C2}\sqrt{\sum_{F\in \cF} \ell(F)\|    [ \QtoG ]_F\times \nu_F \|_{L^2(F)}^2}\right)^2\notag
\end{align}
as a GUB $\|\nabla u - \QtoG\| \le \eta(\T, \QtoG)$ under minimal assumptions on the approximation $\QtoG\in H^1(\T;\R^n)\subset L^2(\Omega; \R^n)$. The constants $C_1$, $C_2$, and $C_H$ (or upper bounds thereof) are computable; cf. Table \ref{tab:expl_const} for an example in 2D with details in Example \ref{ex:constants} at the end of Section \ref{sec:blocks}.
The first assumption is a weakened discrete solution property \begin{equation}\label{eq:solution_property}
    (\QtoG, \nabla w_C)_{L^2(\Omega)}=(f,w_C)_{L^2(\Omega)}\quad\text{for all }w_C\in S_0^1(\cT).
\end{equation}
The second assumption is the orthogonality to the lowest-order divergence-free Raviart-Thomas functions
\begin{align}\label{eq:div_0_property}
    (\QtoG, r)_{L^2(\Omega)}=0\quad\text{for all } r\in RT_0(\T)\cap H(\ddiv=0, \Omega).
\end{align}
\begin{theorem}[residual-based GUB]
\label{thm:reliability}
Suppose that $\QtoG\in H^1(\cT;\R^n)$ and $f\in L^2(\Omega)$ satisfy \eqref{eq:solution_property}--\eqref{eq:div_0_property}.
Then the error estimator
$\eta(\T, \QtoG)$ from \eqref{eq:eta}
is a GUB
$$
\|\nabla u-\QtoG\| \le \eta(\T, \QtoG)
$$
of the error $\|\nabla u-\QtoG\|$ for the solution $u\in V$ to \eqref{eq:poisson}.
The constants $\Ca, C_{2}, \CH$ exclusively 
depend on $\Omega$ and the shape-regularity of $\cT$.
\end{theorem}

The remaining parts of this section are devoted to the proof of Theorem~\ref{thm:reliability} and the computation of (upper bounds of) the constants $\Ca, \Cb$, and $\CH$ in \eqref{eq:eta}. %
The point of departure is the subsequent decomposition that appears necessary in the 
nonconforming and mixed finite element a~posteriori error analysis.
It leads to a split of the error $\|\nabla u-\QtoG\|$ into some divergence part and some consistency part. 

\begin{lemma}[decomposition]\label{lemmadecomposition}
Any $v\in V$ and $\QtoG\in L^2(\Omega; \R^n)$ satisfy the decomposition
\begin{align}\label{eqn:decomp}
\|\nabla v-\QtoG \|^2 = \trb{v-w}^2+\| \QtoG- \nabla w\|^2
\end{align}
with the (unique) minimizer $w\in V$ of the distance 
\[
\delta:= \min_{\varphi \in V}\|\QtoG- \nabla \varphi\|
\] 
of $\QtoG$ to the gradients $\nabla V$ of Sobolev functions.
The solution $u\in V$ to \eqref{eq:poisson} satisfies
\begin{align}
\mu\coloneqq&\;\trb{ f+\ddiv \QtoG}_{*}=\trb{u-w}&&\text{and} \notag\\\label{eqn:decomposition_full}
\|\nabla u-\QtoG \|
^2 =&\;\trb{ f+\ddiv \QtoG}_{*}^2+\| \QtoG- \nabla w\|
^2=\mu^2+\delta^2.
\end{align}
\end{lemma}
\begin{proof}
The minimizer $w \in V$ of $\|\QtoG - \nabla \varphi\|$ among $\varphi \in V$ satisfies the variational formulation $a(w, \varphi)_{L^2(\Omega)} = (\QtoG,\nabla \varphi)_{L^2(\Omega)}$ for all $\varphi \in V$.
(Notice that $w$ is the unique weak solution to the Poisson model problem $-\Delta w = -\ddiv \QtoG \in V^*$.)
In particular, $\QtoG- \nabla w$ is $L^2$ orthogonal onto $ \nabla V$ and the Pythagoras
theorem proves \eqref{eqn:decomp}.
Given  $\varphi\in V$ with $ \trb{\varphi} =1$, the orthogonality of $\QtoG-\nabla w$ to $\nabla\varphi$ and \eqref{eq:poisson} show
\begin{align}\label{eqn:split}
a(u-w,\varphi)_{L^2(\Omega)}= (\nabla u- \QtoG,\nabla\varphi)_{L^2(\Omega)}=\langle f+\ddiv \QtoG,\varphi\rangle
\end{align}
with the duality brackets $\langle\bullet, \bullet\rangle$ in $V^* \times V$.
Since the supremum of \eqref{eqn:split} over all $\varphi\in V$ with $ \trb{\varphi} = 1$ is equal to  $\trb{u-w}= \trb{ f+\ddiv \QtoG}_{*} $,
this and \eqref{eqn:decomp} conclude the proof of \eqref{eqn:decomposition_full}.
\qed
\end{proof}
The split \eqref{eqn:decomp} of the error $\|\nabla u-\QtoG \|$ allows for and enforces a separate estimation of the equilibrium and consistency contribution in residual-based a posteriori error estimators.

In order to derive explicit constants, two lemmas are recalled.
The first has a long tradition in the a~posteriori error control in form of a Helmholtz
decomposition on simply connected domains \cite{c43f5cd9,alonso1996error} and introduces the constant $\CH$ from Theorem \ref{thm:reliability}. 
The following version includes the general case of multiply
connected domains as in \cite{GirRav:86} for $n=2$ or $n=3$ dimensions and weak assumptions 
on a divergence-free function 
$\varrho\in H(\ddiv=0,\Omega)$. 
\begin{lemma}[Helmholtz-decomposition]\label{lem:lemmacc2a}
Suppose the divergence-free function $\varrho\in H(\ddiv=0,\Omega)$ is $L^2$ orthogonal onto $RT_0(\cT)\cap H(\ddiv=0, \Omega)$.
Then there exists $\beta\in H^1(\Omega;\R^N)$, $N=2n-3$, such that any $\beta_C\in S^1(\cT)^N$
satisfies 
\begin{align}\label{eqn:lemmacc2a}
\| \varrho\|^2=\int_\Omega \varrho\cdot\mathrm{Curl}(\beta-\beta_C)\ \mathrm{d} x
\quad\text{and}\quad
\trb{\beta}\le \CH\, \| \varrho\|.
\end{align}
The constant $\CH>0$ exclusively depends on $\Omega$. 
\end{lemma}
\begin{proof}
The compact polyhedral boundary $\partial\Omega$ of the bounded Lipschitz domain $\Omega$
has $J+1$ connectivity components $\Gamma_0,\dots,\Gamma_J$ for some finite
$J\in\bN_0$. Those connectivity components have a positive surface measure $|\Gamma_j|$ and 
a positive distance of each other. So the 
integral mean 
\[
\gamma_j:= \int_{\Gamma_j} \varrho\cdot\nu\, \ \mathrm{d} s /|\Gamma_j|
\]
is well defined and depends continuously on $\varrho\in H(\ddiv=0,\Omega)$ in the sense that
$|\gamma_j|\le c_1\| \varrho\|$ (recall $\ddiv\varrho=0$) for each $j=0,\dots, J$ and $c_1>0$. This constant $c_1$ and the constants $c_2,c_3,c_4$ below exclusively depend on the domain $\Omega$.
The finite real numbers 
$\gamma_0,\dots,\gamma_J$ define the Neumann data for the harmonic 
function $z\in H^1(\Omega)/\R$ with  
\[
\Delta z=0 \text{ in }\Omega \quad\text{and}\quad
\partial z/\partial\nu = \gamma_j\text{ on }\Gamma_j \text{ for all }j=0,\dots, J.
\]
The elliptic regularity theory for polyhedral domains lead to $z\in H^{1+\alpha}(\Omega)$
for some $\alpha>1/2$ and $c_2>0$ with $\|  z \|_{H^{1+\alpha}(\Omega)}\le c_2\, (|\gamma_0|+\dots +|\gamma_J|)$.
The Raviart-Thomas interpolation operator defines a bounded linear operator on
$H(\ddiv,\Omega)\cap L^p(\Omega;\R^n)$ for $p>2$.
It is generally accepted that, for $\alpha>0$ and 
$ \nabla  z\in H(\ddiv=0,\Omega)\cap H^\alpha(\Omega;\R^n)$, the Fortin interpolation 
$I_{\rm F}  \nabla  z \in RT_0(\cT)\cap H(\ddiv=0,\Omega)$ is well defined 
and $\| I_{\rm F}  \nabla  z\|\le c_3 \|  \nabla  z\|_{H^\alpha(\Omega)}$ follows for some $c_3>0$.
The additional property  
$ \nabla  z\in L^p(\Omega;\R^n)$ for some $p>2$ allows the definition of
$\int_F  \nabla  z\cdot \nu_F  \ \mathrm{d} s$ as a Lebesgue integral over a facet $F\in\F$.
One consequence for the boundary facets is the vanishing integral
\[
\int_{\Gamma_j} (\varrho-I_{\rm F}  \nabla  z)\cdot\nu\, \mathrm{d} s =0\quad\text{for all }j=0,\dots, J.
\]
Since $\varrho-I_{\rm F}  \nabla  z\in H(\ddiv=0,\Omega)$ is divergence-free, Theorems 3.1 and 3.4 in \cite{GirRav:86} prove the existence of $c_4>0$ and $\beta\in H^1(\Omega;\R^N)$
with 
\[
 \varrho=I_{\rm F}  \nabla  z+  \Curl\beta
\quad\text{and}\quad
\trb{\beta}\le c_4 \|  \varrho-I_{\rm F}  \nabla  z\|.
\] 
Recall that $\varrho\perp I_{\rm F}\nabla z$ and
$\varrho\perp\curl\beta_C\in RT_0(\T)\cap H(\ddiv=0, \Omega)$. This concludes the proof of
\eqref{eqn:lemmacc2a}
with $\CH\coloneqq\sqrt{1 + (c_1c_2c_3(1+J))^2}c_4$.\qed
\end{proof}

The subsequent version of the trace inequality on the facets $\F$ leads to the piecewise constant $\ell\in P_0(\F)$ defined by
	$$\ell(F)\coloneqq \begin{cases}{}
		(n+1)h_{T}^{2}|F|/|T|&\text{for }F\in\F(\partial\Omega)\cap \F(T),\\
		(n+1) |F|/(h_{T_+}^{-2}|T_+| + h_{T_-}^{-2}|T_-|)&\text{for }F=\partial T_+\cap \partial T_-\in\F(\Omega).
	\end{cases}$$
\begin{lemma}[trace inequality]\label{lem:Skeletal trace inequality}
    Any $f\in H^1(\Omega)$ satisfies
	$$\sum^{}_{F\in\F} \ell(F)^{-1}\left\Vert f\right\Vert_{L^2(F)}^{2} \leq \left\Vert
	h_\T^{-1}f\right\Vert
	^{2} + \frac{2\Ctr}{n}\left\Vert
	h_\T^{-1}f\right\Vert
	\ENorm{f}$$
 with the constant $\Ctr\coloneqq \max_{T\in\T}\max_{x\in T}|x-\mathrm{mid}(T)|/h_T<n/(n+1)$.

\end{lemma}
\begin{proof}
	The center of inertia $\mathrm{mid}(T)=\sum_{j=0}^n P_j/(n+1)$ of the $n$-simplex $T=\mathrm{conv}\{P_0, ..., P_n\}\in\T$ and the $n+1$ faces $F_j=\mathrm{conv}\{P_0, ..., P_{j-1},P_{j+1}, ..., P_n\}\in \mathcal F(T)$ for $j=0, ..., n$ give rise to the decomposition of $T$ into $n+1$ sub-simplices $T_j' = \mathrm{conv}(F_j, \mathrm{mid}(T))$ with volume 
	$|T_j'|=|T|/(n+1)$.

	Standard arguments like the trace identity on $T_j'\subset T$ \cite[Lemma
	2.1]{carstensen_explicit_2012} for $|f|^2\in
	W^{1,1}(T')$ and a Cauchy inequality show 
	$$\frac1{|F_j|}\left\Vert f\right\Vert_{L^2(F_j)}^{2}\leq \frac{1}{|T_j'|}\left\Vert
	f\right\Vert_{L^2(T_j')}^{2}+ \frac{2}{n|T_j'|}\max_{x\in T_j'}|x-\mathrm{mid}(T)|\left\Vert
	f\right\Vert_{L^2({T_j'})}^{}\|\nabla f\|_{L^2(T')}.$$
	The distance $\max_{x\in T_j'}|x-\mathrm{mid}(T)|=|P_k - \mathrm{mid}(T)|\leq \Ctr h_T$ is attained at a vertex $P_k$ for $k\in \{0, ..., j-1, j+1, ..., n\}$.
	Since the centroid $\mathrm{mid}(T)$ divides each median of $T$ in the ratio $n$ to $1$ and the length of each median
	is strictly bounded by $h_T$, 
	the bound $\Ctr<n/(n+1)$ follows and cannot be improved in the absence of further assumptions on the shape of the simplex $T$.
	Since $|T_j'| = |T|/(n+1)$, %
	the previously displayed estimate leads to
	\begin{align*}
		\frac{|T|}{(n+1)h_{T}^2|F_j|}\left\Vert f\right\Vert_{L^2(F_j)}^{2}&\leq\left\Vert
	h_{T}^{-1}f\right\Vert_{L^2(T_j')}^{2}+ \frac{2\Ctr}{n}\left\Vert
	h_T^{-1}f\right\Vert_{L^2({T_j'})}^{}\|\nabla f\|_{L^2(T_j')}.%
	\end{align*}
	Let $\T'$ be the refinement of $\T$, obtained by replacing $T\in\T$ with $T_0', ..., T_d'$ from above.
	The triangulation $\T'$ allows for the facet based decomposition $\{\omega'(F)\}_{F\in\F}$ of $\Omega$, where $\omega'(F)$
	is either the patch $\omega'(F)=\mathrm{int}(T'_+\cup T'_-)$ for an interior facet
	$F=T'_+\cap T'_-$ or $\omega'(F)=\mathrm{int}(T')$ for
	$F\in\F(\partial\Omega)\cap\F(T')$.
    This establishes, for any $F\in\F$, the estimate 
	\begin{align*}
		\ell(F)^{-1}\left\Vert f\right\Vert_{L^2(F)}^{2}
		&\leq \left\Vert
		h_{\T}^{-1}f\right\Vert_{L^2(\omega'(F))}^{2} + \frac{2\Ctr}{n}\left\Vert
		h_{\T}^{-1}f\right\Vert_{L^2(\omega'(F))}\|\nabla f\|_{L^2(\omega'(F))}.
	\end{align*}
	Since the family $\{\omega'(F) : F\in\F\}$ has no overlap, the sum of the last displayed inequality  over all $F \in \F$ and a Cauchy inequality conclude the proof of Lemma \ref{lem:Skeletal trace inequality}.\qed
\end{proof}
The next lemma utilizes a quasi-interpolation operator $J: H^1(\Omega) \to S^1(\T)$ with the restriction $J(V) \subset S^1_0(\T)$, e.g., $J = J_1 \circ I_{NC}$ from \cite[Section 5]{carstensen_constants_2018} with explicit constants for $n=2$, and the approximation and stability properties 
\begin{align}\label{eqn:approx_stab}
    \left\Vert h_\T^{-1}(\varphi-J\varphi)\right\Vert
    &\leq \cnst{cnst:C1} \ENorm{\varphi}\quad\text{and}\quad
    \ENorm{\varphi-J\varphi}\leq \cstb \ENorm{\varphi}
\end{align}
for constants $\cnst{cnst:C1}$ and $\cstb$ exclusively depending on the shape-regularity of $\T$.
For the precise definition of $J_1$ and $I_{NC}$, we refer to \cite[eq.~(47) and Section 5]{carstensen_constants_2018}.
Recall the constant $\Ctr$ from Lemma \ref{lem:Skeletal trace inequality} and set $\cnst{cnst:C2}\coloneqq (\cnst{cnst:C1}(\cnst{cnst:C1}+2\Ctr\cstb\,/n))^{1/2}$.
\begin{lemma}[equilibrium]\label{lem:eq1} 
Suppose that $\QtoG \in H(\ddiv,\cT )$ and $f\in L^2(\Omega)$ satisfy \eqref{eq:solution_property} and suppose
$(\QtoG|_T)|_F\cdot \nu_F\in L^2(F)$ for all $F\in\cF(T)$ and $T\in \cT$.
Then
\[
	\trb{f+\ddiv \QtoG}_*  \le \cnst{cnst:C1}\|  h_\cT( f+{\ddiv}_{\rm pw} \QtoG) \|
+ \cnst{cnst:C2}\sqrt{\sum_{F\in\cF(\Omega)} \ell(F) \|  [ \QtoG ]_F\cdot \nu_F  \|_{L^2(F)}^2} .
\]
\end{lemma}

\begin{proof}
Given  $\varphi\in V$ with $ \trb{\varphi} =1$, set $\psi:=\varphi-\varphi_C$ for some quasi-interpolation $\varphi_C\coloneqq J\varphi\in S^1_0(\cT)$ with \eqref{eqn:approx_stab}.
Since \eqref{eq:solution_property} implies $\langle f+\ddiv \QtoG,\varphi\rangle =\langle f+\ddiv \QtoG,\psi\rangle$, a piecewise integration by parts and the collection of jump contributions show
\begin{align}
    \langle f+\ddiv \QtoG,\varphi\rangle =
    (f + \divpw \QtoG, \psi)_{L^2(\Omega)}
    -\sum_{F\in\cF(\Omega)} \langle[\QtoG]_F\cdot \nu_F, \psi\rangle_{L^2(F)}.
    \label{eq:proof-equilibrium-ibp}
\end{align}
The first bound follows from a Cauchy inequality and \eqref{eqn:approx_stab},
\begin{align}
    (f + \divpw \QtoG, \psi)_{L^2(\Omega)} &\leq \|h_\T(f+\divpw \QtoG)\|\, \|h_\T^{-1}\psi\|
    \label{ineq:proof-equilibrium-1}\\
    &\leq \cnst{cnst:C1} \|h_\T(f+\divpw \QtoG)\|\,
	\trb{\varphi}.\nonumber
\end{align}
The second bound additionally exploits the trace inequality of
Lemma \ref{lem:Skeletal trace inequality},
\begin{align}
    \sum_{F\in\cF(\Omega)} &\langle[\QtoG]_F\cdot \nu_F, \psi\rangle_{L^2(F)}
    \label{ineq:proof-equilibrium-2}\\
    &\leq \sqrt{\sum_{F\in\cF(\Omega)} \ell(F)\|[\QtoG]_F\cdot \nu_F\|_{L^2(F)}^2} \sqrt{\sum_{F\in \cF(\Omega)}\ell(F)^{-1}\|\psi\|_{L^2(F)}^2}\nonumber\\
															  &\leq \sqrt{\cnst{cnst:C1}^2 +
															  \frac{2\Ctr}{n}\cstb\cnst{cnst:C1}}\sqrt{\sum_{F\in\cF(\Omega)}
																  \ell(F)\|[\QtoG]_F\cdot \nu_F\|_{L^2(F)}^2}\nonumber
\end{align}
with $\ENorm{\varphi}=1$ in the last step.
Since \eqref{ineq:proof-equilibrium-1}--\eqref{ineq:proof-equilibrium-2} hold for all $\varphi\in V$ with $\trb{\varphi}=1$, the supremum in \eqref{eq:proof-equilibrium-ibp} over all such $\varphi$ concludes the proof.\qed
\end{proof}
The final ingredient for the proof of Theorem \ref{thm:reliability} controls the second term $\delta$ in the decomposition of Lemma \ref{lemmadecomposition} for $\varrho\coloneqq \QtoG-\nabla w$. Recall $\CH$ from Lemma \ref{lem:lemmacc2a} and $\cnst{cnst:C1}$ from \eqref{eqn:approx_stab}, and $\cnst{cnst:C2}$ from page 9.

\begin{lemma}[conformity]\label{lem:eq2}
Suppose the divergence-free function $\varrho\in H(\ddiv=0,\Omega)\cap H(\curl,\cT)$ is $L^2$ orthogonal onto $RT_0(\cT) \cap H(\ddiv=0,\Omega)$ and satisfies 
$(\varrho|_T)|_F\times \nu_F\in L^2(F)$ for all $F\in\cF(T)$ and $T\in T$.
Then
\[
\CH^{-1}\|\varrho\| \le \cnst{cnst:C1}\|  h_\cT{\curl}_{\rm pw} \varrho \|
+  \cnst{cnst:C2}\sqrt{\sum_{F\in\cF} \ell(F) \|  [  \varrho]_F \times \nu_F \|_{L^2(F)}^2}.\]
\end{lemma}

\begin{proof}
Lemma \ref{lem:lemmacc2a} provides $\beta\in H^1(\Omega; \R^N)$ with 
\eqref{eqn:lemmacc2a}
for a 
(component-wise) quasi-interpolation $\beta_C\in S^1(\T)^N$ with \eqref{eqn:approx_stab} as in the proof of Lemma \ref{lem:eq1}; set $\psi:=\beta-\beta_C$.
A piecewise integration by parts and the collection of jump contributions shows
\[
 \|  \varrho\| ^2 = \int_\Omega \varrho\cdot\Curl(\beta-\beta_C)dx =
\int_\Omega \psi \cdot {\curl}_{\rm pw} \varrho \, dx 
+\sum_{F\in\cF} \int_F  \psi [  \varrho]_F\times\nu_F \, ds.
\]
Stability and approximation properties of the  
quasi-interpolation \eqref{eqn:approx_stab} and the trace inequality of Lemma \ref{lem:Skeletal trace inequality} eventually lead to
\begin{align*}
	\left\Vert \varrho\right\Vert^{2}& \leq \cnst{cnst:C1}\left\Vert h_\T{\curl}_{\rm pw} \varrho\right\Vert_{}^{}\trb{\beta} + \cnst{cnst:C2} \sqrt{\sum_{F\in\cF}
\ell(F)\|[\varrho]_F\times \nu_F\|_{L^2(F)}^2}\trb{\beta}.
\end{align*}
In fact, the routine estimation with element and jump terms is completely analogous to the proof of Lemma \ref{lem:eq1} and leads to the same constants $\cnst{cnst:C1}, \cnst{cnst:C2}$.
This and $\trb{\beta}\leq \CH\|\varrho\|$ conclude the proof.\qed
\end{proof}

\begin{proof}[Theorem \ref{thm:reliability}]
The trace of $\QtoG|_T\in H^1(T; \R^n)$ is well defined on any facet $F\in\F(T)$ of the simplex $T\in\T$. 
Lemma \ref{lemmadecomposition} provides $w\in V=H^1_0(\Omega)$ with $\|   \nabla u-\QtoG \|^2 = \trb{ f + \ddiv \QtoG}_{*}^2+\| \QtoG- \nabla w\|^2$. 
Since $\QtoG$ satisfies \eqref{eq:solution_property}, Lemma \ref{lem:eq1} establishes
$$\trb{ f + \ddiv \QtoG}_{*}  \le \CeqT\|  h_\cT( f+{\ddiv}_{\rm pw} \QtoG) \|
+ \CeqF\sqrt{\sum_{F\in\cF(\Omega)} \ell(F) \|  [ \QtoG ]_F\cdot \nu_F  \|_{L^2(F)}^2}.$$
The assumption \eqref{eq:div_0_property} on $\QtoG$ and an integration by parts prove that Lemma \ref{lem:eq2} is applicable to $\varrho\coloneqq \QtoG-\nabla w\in H(\ddiv=0, \Omega)\cap H(\curl, \T; \R^n)$.
Since $\curl\nabla w = 0$ and $\nabla w\times \nu = 0$, 
this reveals 
$$\CH^{-1}\| \QtoG- \nabla w\|\leq\CeqT\|  h_\cT{\curl}_{\rm pw} \QtoG \|
+  \CeqF\sqrt{\sum_{F\in\cF} \ell(F) \|  [\QtoG]_F \times \nu_F \|_{L^2(F)}^2}.
$$
The above estimates together with the decomposition of Lemma \ref{lemmadecomposition} establish $\eta(\T, \QtoG)$ as a GUB for the error $\|\nabla u - \QtoG\|$. \qed
\end{proof}

\begin{example}[constants for right-isosceles triangles]\label{ex:constants}
In two space dimensions, $\|\Curl \bullet\| = \trb{\bullet}$ and so $\CH\leq 1$ for a simply connected domain $\Omega$ in Lemma \ref{lem:lemmacc2a}.
The choice $J\coloneqq J_1\circ I_{\rm NC}$ from \cite[Section 5]{carstensen_constants_2018} of the quasi-interpolation operator $J$ in the proof of Lemma \ref{lem:eq1} allows for the explicit estimates $$\cnst{cnst:C1}\leq\sqrt{48^{-1} + \jj^{-2} + \capx^2} \text{ and } \cstb\leq 1+\sqrt{72}\capx,$$
where $j_{1,1} = 3.8317$ denotes the first positive root of the first Bessel function.
For
triangulations into right-isosceles triangles, the constant $\capx\leq\sqrt{3}/(2 - 2\cos(\pi/\max\{4,M_{\rm bd}\}))$ %
from \cite[Lemma 4.8]{carstensen_constants_2018} depends on the domain by the maximal number $M_{\rm bd}\leq 4\max\{\pi,\omega_{\rm max}\}/\pi\leq 8$ of triangles sharing a boundary vertex. %
Given the maximal interior angle $\omega_{\rm max}$ of $\Omega$, Table \ref{tab:expl_const} displays those constants for the maximal possible value $M_{\rm bd}= 4\max\{\pi,\omega_{\rm max}\}/\pi$.
The geometric quantity $\max_{x\in T}|x-\mathrm{mid}(T)|$ equals two-thirds of the maximum median of $T$.
Thus, $\Ctr=\sqrt{5}/(3\sqrt{2})\leq 0.5271$ and $\ell(F) = 6h_F$ for interior edges $F\in\F(\Omega)$ and $\ell(F)=12h_F$ for boundary edges $F\in\F(\partial\Omega)$ of triangulations into right-isosceles triangles.
Consequently,
\begin{align}\label{eqn:C_T}
    \cnst{cnst:C1}&\leq \sqrt{48^{-1} + \jj^{-2} + \capx^2}=:C_\T,\\
    \cnst{cnst:C2}&\leq \sqrt{C_\T(C_\T+0.5271(1+\sqrt{72}\capx))}=:C_\F\label{eqn:C_F}.
\end{align}

\end{example}

\section{Explicit residual-based a posteriori HHO error estimator}
\label{sec:residual}
The arguments from Section \ref{sec:blocks} apply to the HHO method and result in a stabilization-free reliable a posteriori error control. 
In combination with the efficiency estimate from this section, this leads to a new explicit residual-based a posteriori error estimator for the HHO method that is equivalent to the error up to data oscillations.
\subsection{Hybrid high-order methodology}\label{sec:HHO}
The HHO ansatz space reads $\Vh \coloneqq P_k(\T) \times P_k(\F(\Omega))$
for $k\in\N_0$ with the subspace $P_k(\F(\Omega)) \subset P_k(\F)$ of piecewise polynomials $p\in P_k(\F)$ under the convention $p_{|\partial\Omega}=0$.
The interpolation $\mathrm{I}: V \to V_h$ maps $v \in V$ onto $\mathrm{I} v \coloneqq (\Pi_k v, \Pi_{\mathcal{F},k} v) \in V_h$.
Given any $v_h = \vTF \in \Vh$, the reconstruction operator $R:\Vh \rightarrow P_{k+1}(\T)$ defines the unique piecewise polynomial $Rv_h \in P_{k+1}(\T)$ with $\Pi_{0}(Rv_h - \vT) = 0$ such that, for all $w_{k+1} \in P_{k+1}(\T)$,
\begin{align}
\label{eqn:R}
    &\apw(Rv_h,w_{k+1})\nonumber\\
    &\qquad= \apw(\vT, w_{k+1})
    -\sum_{T\in\T}\langle v_{\T|T}-\vF, \nabla w_{k+1|T} \cdot \nu_T \rangle_{L^2(\partial T)}.
\end{align}
Let $u_h\in \Vh$ solve the HHO discrete formulation of \eqref{eq:poisson} with
\begin{align}\label{eqn:HHO}
a_h(u_h,v_h)= (f,\vT)_{L^2(\Omega)} \qquad \text{for all } v_h=\vTF \in \Vh
\end{align}
for the HHO bilinear form 
\begin{align}%
    a_h(u_h,v_h) &\coloneqq\apw( Ru_h,  Rv_h)+
 s_h(u_h, v_h)%
\end{align}
and the stabilization term $s_h(u_h, v_h)$ from \eqref{eqn:s_h0}.
Given any $w_C \in S^1_0(\T)= P_1(\T)\cap H^1_0(\Omega)$,
the definition of the reconstruction operator $R$ in \eqref{eqn:R} verifies $R \mathrm{I} w_C=w_C$ with the interpolation $\mathrm{I}$ onto $V_h$.
Hence, $\SpTF \mathrm{I} w_C=0$ vanishes for all $F \in \F(T)$ and $T \in \T$.
This and \eqref{eqn:HHO} show, for all $w_C\in S^1_0(\T)$, that
\begin{align}\label{eqn:HHO_solution_property}
    \apw(Ru_h, w_C) = (f, \Pi_k w_C)_{L^2(\Omega)} = (\Pi_k f, w_C)_{L^2(\Omega)}.
\end{align}
\subsection{Explicit a posteriori error estimator}\label{sub:HHO_eta}
As a result of \eqref{eqn:HHO_solution_property}, $\nablapw Ru_h$ satisfies the solution property \eqref{eq:solution_property} if $k\geq 1$ and, in the lowest order case $k=0$, \eqref{eq:solution_property} holds with $f$ replaced by $\Pi_0 f$.
This allows the application of the theory from Section \ref{sec:blocks} to the HHO method with minor modifications for the case $k=0$.
Define the error estimator contributions
\begin{equation}\label{eqn:etares_parts}
    \begin{split}
         \eta^2_{\text{res}, 1}(\T)&\coloneqq \begin{cases}\|h_\T(  f+\Delta_{\rm pw} Ru_h)\|^2&\text{for } k\geq1,\\
    \|h_\T\Pi_0  f\|^2&\text{for } k=0,
    \end{cases}\\
        \eta^2_{\text{res}, 2}(\T)&\coloneqq \begin{cases}0&\text{for } k\geq1,\\
    \osc0^2(f, \T)&\text{for } k=0,
    \end{cases}\\
	\eta^2_{\text{res}, 3}(\T)&\coloneqq\sum_{F\in \cF(\Omega)} \ell(F)\|    [ \nablapw Ru_h ]_F\cdot \nu_F \|_{L^2(F)}^2,\\
    \eta^2_{\text{res}, 4}(\T)&\coloneqq\sum_{F\in \cF} \ell(F)\|    [ \nablapw Ru_h ]_F\times \nu_F \|_{L^2(F)}^2.
    \end{split}
\end{equation}
Since $\nablapw Ru_h$ is a piecewise gradient, its piecewise $\curl$ vanishes.
This leads to the explicit residual-based a posteriori error estimator 
\begin{equation}
\label{eq:residualdef}
\eta_{\text{res}}^2(\T)\coloneqq \left(\cnst{cnst:C1}\eta_{\text{res}, 1}(\T) + C_P\eta_{\text{res}, 2}(\T)
+\cnst{cnst:C2}\eta_{\text{res}, 3}(\T)\right)^2+ \CH^2\cnst{cnst:C2}^2\eta_{\text{res}, 4}^2(\T).
\end{equation}
(Recall $\cnst{cnst:C1}, \cnst{cnst:C2} $ from Lemma \ref{lem:eq1} and $\CH$ from Lemma \ref{lem:lemmacc2a} as well as the Poincar\'e constant $C_P\leq \pi^{-1}$.)
The main result of this section verifies the assumptions in Theorem \ref{thm:reliability} and proves reliability and efficiency of $\etares(\T)$.
\begin{theorem}[residual-based GUB for HHO]\label{thm:HHO_equivalence}
Let $u\in V$ solve the 
Poisson equation \eqref{eq:poisson} and let
$u_h \in V_h$ solve the discrete 
formulation \eqref{eqn:HHO}.
Then
\begin{align*}
	\npw{u-\cR u_h} \leq \etares(\T)\leq \Ceff\big( \npw{u-\cR u_h}+\osc q(f, \cT)\big)
\end{align*}
and $\osc {k-1}(f, \cT)\leq \cnst{cnst:oscres_rel} \etares(\T)$ hold for any $q\in\N_0$.
The constants $\newcnst\label{cnst:res_eff}$ and $\newcnst\label{cnst:oscres_rel}$ exclusively depend on $k, q$ and on the shape-regularity of the triangulation $\cT$.
\end{theorem}
\subsection{Proof of Theorem \ref{thm:HHO_equivalence}}
The orthogonality of $\nabla_\pw R u_h$ to the divergence-free Raviart-Thomas space {of} lowest degree is an assumption in Theorem \ref{thm:reliability} and verified below.
\begin{lemma}[orthogonality]
\label{lem:GuhorthRTd0} 
The piecewise gradients $\nablapw R\Vh$
are $L^2$ orthogonal to the space {~$RT_0(\T)\cap H(\mathrm{div} = 0,\Omega)$}, i.e., any $v_h\in\Vh$ and $\qRT \in RT_0(\T)\cap H(\mathrm{div} = 0,\Omega)$ satisfy
\begin{align}
(\nablapw Rv_h, \qRT)_{L^2(\Omega)} = 0.
\end{align}
\end{lemma}
\begin{proof}
Given any $\qRT \in RT_0(\T)\cap H(\mathrm{div} = 0,\Omega)$, $\ddiv \qRT=0$ shows
 $\qRT \in P_0(\T;\R^n)$ \cite[Lemma 14.9]{ErnGuermond2021}. 
 Since $P_0(\T;\R^n) = \nablapw P_1(\T)$, there exists a piecewise affine function $\phi_1 \in P_1(\T)$ with $\qRT =  \nabla_\pw  \phi_1$ a.e.~in $\Omega$.
This and the definition of $\cR  v_h$ from \eqref{eqn:R} imply, for any $v_h = \vTF\in V_h$, that
\begin{align*}
(\nablapw Rv_h,\qRT)_{L^2(\Omega)} &= a_\pw(Rv_h, \phi_1)\\
&=  a_\pw(v_\T,\phi_1) - \sum\limits_{T\in \mathcal T} \langle  \vT|_T-\vF, \nabla_\pw  \phi_1 \cdot \bn_T \rangle_{L^2(\partial T)}.
\end{align*}
This, a piecewise integration by parts, and $\Delta_\pw \phi_1 \equiv 0$ lead to
\begin{align}
(\nablapw \cR v_h,\qRT)_{L^2(\Omega)} &= \sum\limits_{T\in \mathcal T}\langle  \vF,\nabla_\pw  \phi_1 \cdot \bn_T \rangle_{L^2(\partial T)}\nonumber\\
&= \sum\limits_{F\in \mathcal F}\langle  \vF, [\qRT\cdot\bn_F]_F \rangle_{L^2(F)}.\label{eq:proof-orthogonality-RT-div=0-ibp}
\end{align}
Since $\qRT\in RT_0(\T)$ has continuous normal components, the jump term $\langle \vF, [\qRT]_F\cdot \bn_F\rangle_{L^2(F)} =0$ vanishes for all $F\in \facetsI$. This, $\vF\equiv 0$ on $\partial\Omega$, and \eqref{eq:proof-orthogonality-RT-div=0-ibp} conclude
$(\nablapw \cR v_h,\qRT)_{L^2(\Omega)} = 0$. \qed
\end{proof}
The following lemma concerns the efficiency of the jump contributions. %
Each facet $F\in\F$ has at most two adjacent simplices that define a triangulation $\T(F)\coloneqq\{T\in\T : F\in\F(T)\}$ of the facet-patch $\omega(F)\coloneqq \mathrm{int}(\bigcup_{T\in\T(F)}T)$.

\begin{lemma}[efficiency of jumps]
\label{lem:efficiency}
The solution $u \in V$ to \eqref{eq:poisson} and the discrete solution $u_h \in V_h$ to \eqref{eqn:HHO} satisfy (a) for all $F \in \F$ and (b) for all $F \in \F(\Omega)$.
\begin{enumerate}
    \item[(a)] $h_F^{1/2}\|[\nablapw\cR u_h]_F\times \nu_F \|_{L^2(F)} \lesssim \min_{v \in V} \|\nabla v -  \nablapw \cR u_h\|_{L^2(\oF)}$,\\
    \item[(b)] $h_F^{1/2}\|[\nablapw\cR u_h]_F\cdot \nu_F \|_{L^2(F)} \lesssim \|\nabla  u -  \nablapw \cR u_h\|_{L^2(\oF)} + \osck(f, \T(F))$.
\end{enumerate}
\end{lemma}
\begin{proof}
\newcommand{\conv}{\mathrm{conv}}
The proof is based on the following extension argument.
Given a polynomial $p\in P_k(F)$ of degree at most $k$ along the side $F\in\F$, the coefficients determine a polynomial (also denoted by $p$)  
along the hyperplane $H$ that enlarges $F$. The intersection $\widehat{F}:= H\cap \conv({ \omega(F)}) $  of the hyperplane 
$H$ with the convex hull of the facet-patch $\omega(F)$ may be strictly larger than $F$.
The shape-regularity  of $\T$ bounds the size of  $\widehat{F}$ in terms of $F$ and an inverse estimate leads to a bound $\| p \|_{L^\infty(\widehat{F})} \le C(k)\| p \|_{L^\infty(F)}$
with a constant $C(k)$ that depends on the shape-regularity of $\T$ and on $k$. The extension of $p$ from $H$ to $\R^n$ by constant values along the side normal $\nu_F$ leads to a polynomial $\widehat{p}\in P_k(\R^
n)$ with 
\begin{align}\label{eqn:L_infty_p}\| \widehat{p} \|_{L^\infty(\omega(F))}
\le
\| p \|_{L^\infty(\widehat{F})} \le C(k)\| p \|_{L^\infty(F)}.
\end{align}
\emph{Proof of (a).}
The %
tangential jump $\varrho_F:= [\nablapw\cR v_h]_F\times \nu_F \in P_{k}(F;\mathbb{R}^N)$
is a polynomial in $N=2n-3$ components on $F\in\F$ for $n=2,3$.
Let $p=\varrho_F(j)$ be one of the components of $\varrho_F\in P_k(F)^N$, for $j=1,...,N$, and extend it as explained above to $\hat p\in P_k(\R^n)$ and call this $\widehat\varrho(j)$ in the vector 
 $\widehat{\varrho_F}\in P_k(\R^n;\mathbb{R}^N)$.
The proof involves the piecewise polynomial facet-bubble function $b_F\coloneqq n^n \Pi_{j=1}^n \varphi_j $ for the $n$
nodal basis function $\varphi_1, \dots, \varphi_n \in S^1(\part)$ associated with the vertices of $F$.
An inverse estimate \cite[Proposition 3.37]{Ver:13} shows
\begin{align}\label{eqn:lemma_7_a}
\| \varrho_F \|_{L^2(F)}^2 \lesssim \| b_F^{1/2}\varrho_F \|_{L^2(F)}^2=\langle b_F\varrho_F, \varrho_F\rangle_{L^2(F)}.
\end{align}
 Since $\varrho:= b_F\widehat{\varrho_F}  \in S^{k+n}_0(\T(F);\mathbb{R}^N)$ vanishes on  $\partial \oF\setminus\mathrm{int}(F)$, \eqref{eqn:lemma_7_a} and a piecewise integration by parts show
 \begin{align}
\|\varrho_F \|_{L^2(F)}^2 \lesssim \langle
\varrho,
[\nablapw\cR v_h]_F\times \nu_F
\rangle_{L^2(F)} = (\curl \varrho, \nablapw R v_h )_{L^2(\oF)}.
\label{ineq:bound-tangential-jump-ibp}
\end{align}
This and $(\curl\varrho, \nabla v)_{L^2(\Omega)} = 0$ for any $v\in V$ imply
\begin{align}
\label{eq:lemma2_s2ls}
\|\varrho_F \|_{L^2(F)}^2 &\lesssim (\curl\varrho,\nablapw(R v_h-v))_{L^2(\oF)}\notag\\
&\leq \| \curl\varrho \|_{L^2(\oF)} \| \nablapw(R v_h-v) \|_{L^2(\oF)}.
\end{align}
An inverse estimate, $\|b_F\|_{L^\infty(\omega(F))} = 1$, and \eqref{eqn:L_infty_p} imply
\begin{align}
 \| \curl\varrho \|_{L^2(\oF)} 
 &\lesssim \|\nabla\varrho\|_{L^2(\oF)} 
  \lesssim h_F^{-1+n/2} \| \varrho \|_{L^\infty(\oF)}\nonumber\\
  &\lesssim h_F^{-1+n/2}  h_F^{-(n-1)/2}\| \varrho_F  \|_{L^2(F)}  = h_F^{-1/2}\| \varrho_F  \|_{L^2(F)}\label{ineq:bound-curl-tangential-jump}.
 \end{align}
In combination with \eqref{eq:lemma2_s2ls}, this concludes the proof of (a).
\medskip

\noindent\emph{Proof of (b).}
The efficiency of normal jumps (b) follows from the arguments for conforming FEMs, cf.~\cite[Section 1.4.5]{Ver:13}; further details are omitted. \qed
\end{proof}
The following lemma reveals that the order $k \geq \mathbb{N}_0$ of the oscillations $\osc k(f, T)$ in Lemma \ref{lem:efficiency} (b) can be any natural number.
It is certainly known to the experts but hard to find in the literature.
Recall the convention $\Pi_{-1}\coloneqq0$.
\begin{lemma}[efficiency of lower-order oscillations]
\label{lem:local-data-oscillation}
    Given any simplex $T \in \mathcal{T}$ and parameters $k, q \in \mathbb{N}_0$,
    the solution $u \in V$ to \eqref{eq:poisson} satisfies
	\begin{align}
		&\cnst{cnst:local-data-oscillation}^{-1}\osc{k-1}^2( f,T)\label{ineq:osc-best-order}
		\leq \min_{v_{k+1} \in P_{k+1}(T)}\|\nabla_\pw(u - v_{k+1})\|_{L^2(T)}^2+\osc{q}^2(f,T).
	\end{align}
	The constant $\newcnst\label{cnst:local-data-oscillation}$ exclusively depends on $q$ and the shape of $T$.
\end{lemma}

\begin{proof}
The assertion \eqref{ineq:osc-best-order} is trivial for $q \leq k-1$,
so suppose $k \leq q$.
Any $v_{k+1}\in P_{k+1}(T)$ and $\varrho_T\coloneqq \Pi_q f + \Delta v_{k+1}\in P_q(T)$ satisfy
\begin{align}\label{eqn:osc_split}
    \osc{k-1}^2(f, T) \leq h_T^2\|f + \Delta v_{k+1}\|_{L^2(T)}^2 = \osc{q}^2(f, T) + h_T^2\|\varrho_T\|_{L^2(T)}^2.
\end{align}
Let $b_T\in S^{n+1}_0(T)$ with $0\leq b_T\leq 1=\max b_T$ denote the volume bubble-function on $T\in\T$.
The equivalence of norms in the finite-dimensional space $P_q(T)$ provides 
\begin{align}\label{eqn:bT_equivalence}
    \|b_T^{1/2}\varrho_T\|_{L^2(T)} \leq \|\varrho_T\|_{L^2(T)}\leq \cnst{cnst:bT_efficiency} \|b_T^{1/2}\varrho_T\|_{L^2(T)}.
\end{align}
A more detailed analysis of the mass matrices reveals that the constant $\newcnst\label{cnst:bT_efficiency}$ exclusively depends on the polynomial degree $q$. An integration by parts with $b_T\varrho_T\in S^{q+n+1}_0(T)\subset V$ and the weak formulation \eqref{eq:poisson} result in
	\begin{align*}
	    \|b_T^{1/2}\varrho_T\|_{L^2(T)}^2 &= (\Pi_q f + \Delta v_{k+1}, b_T\varrho_T)_{L^2(T)}\\
	    &= (\nabla(u-v_{k+1}), \nabla (b_T\varrho_T))_{L^2(T)} - (f - \Pi_q f, b_T\varrho_T)_{L^2(T)}.
	\end{align*}
	A Cauchy inequality, the inverse estimate $h_T\|\nabla (b_T\varrho_T)\|_{L^2(T)}\leq \newcnst\label{cnst:bT_inverse}\|\varrho_T\|_{L^2(T)}$
	with a constant $\cnst{cnst:bT_inverse}$ that exclusively depends on $q+n+1$ and the shape of $T$, and \eqref{eqn:bT_equivalence} lead to
	\begin{align}\label{eqn:osc_efficiency_proof}
	    \cnst{cnst:bT_efficiency}^{-2}h_T\|\varrho_T\|_{L^2(T)}\leq \cnst{cnst:bT_inverse}\|\nabla(u-v_{k+1})\|_{L^2(T)} + \osc q(f, T)\big.
	\end{align}
	The combination of \eqref{eqn:osc_split} with \eqref{eqn:osc_efficiency_proof} and a Cauchy inequality conclude the proof of \eqref{ineq:osc-best-order}, e.g., with $\cnst{cnst:local-data-oscillation}= 1 +  \cnst{cnst:bT_efficiency}^4(1 + \cnst{cnst:bT_inverse}^2).$\qed
\end{proof}

\begin{proof}[of Theorem \ref{thm:HHO_equivalence}]
Recall the definition of $\etares(\T)$ for $k\geq 1$ and $k=0$ in \eqref{eqn:etares_parts}. 
Since
$
    \osc{k-1}^2(f, \T)\leq \eta_{\text{res},1}^2(\T) + \eta_{\text{res},2}^2(\T)\lesssim \etares^2(\T)$,
the remaining parts of this proof discuss the reliability and efficiency of $\etares(\T)$.

Lemma \ref{lem:GuhorthRTd0} provides
the orthogonality of $\nablapw Ru_h\in H^1(\T; \R^n)$ to the divergence-free Raviart-Thomas function $RT_0(\T)\cap H(\mathrm{div} = 0,\Omega)$.
This and \eqref{eqn:HHO_solution_property} show that the assumptions 
in Theorem \ref{thm:reliability} hold for $\QtoG \coloneqq \nablapw Ru_h$ and $k\geq 1$, whence the reliability of $\etares(\T)$ follows with a reliability constant 1. 

Minor modifications to the proof of Theorem \ref{thm:reliability} lead to reliability in the case $k=0$.
In fact, the only modifications required concern the upper bound of $\trb{f + \mathrm{div}\,\nablapw R u_h}_* \leq \trb{(1 - \Pi_0) f}_* + \trb{\Pi_0 f + \mathrm{div}\,\nablapw R u_h}_*$.
A piecewise Poincar\'e inequality shows $\trb{(1 - \Pi_0) f}_* \leq C_P\osc{0}(\T,f)$ with the Poincar\'e constant $C_P$($\leq 1/\pi$ for simplices). Lemma \ref{lem:eq1} proves $\trb{\Pi_0 f + \mathrm{div}\,\nablapw R u_h}_* \leq \Ca\eta_{\text{res},1}(\T) + \Cb \eta_{\text{res},3}(\T)$. Hence, the decomposition of Lemma \ref{lemmadecomposition} and Lemma \ref{lem:eq2} result in $\trb{u-Ru_h}_\pw\leq\etares(\T)$.
This provides the reliability and it remains to verify the efficiency $\etares(\T)\lesssim \trb{u-Ru_h}_\pw + \osc q(f,\T)$ for any $q\in\mathbb N_0$.
The Pythagoras theorem and \eqref{eqn:osc_efficiency_proof} with $\varrho_T\coloneqq \Pi_kf + \Delta R u_h\in P_k(T)$ and $v_{k+1}\coloneqq Ru_h\in P_{k+1}(T)$ lead to the local efficiency of the volume contributions
\begin{align*}
    \|h_T(f+\Delta_\pw Ru_h)\|_{L^2(T)}^2&=\osck^2(f, T) + \|h_T\varrho_T\|_{L^2(T)}^2\\
    &\lesssim \|\nabla(u- Ru_h)\|_{L^2(T)}^2 + \osck^2(f, T).
\end{align*}
Lemma \ref{lem:efficiency} considers the remaining terms in the error estimator and establishes their efficiency namely,
\begin{align*}
    \sum_{F\in\F}h_F\| [ \nabla_\pw R u_h ]_F\|_{L^2(F)}\lesssim \trb{u-Ru_h}_\pw + \osck(f,\T)
\end{align*}
with the modified jump $[\nabla_\pw R u_h]_F = \nabla_\pw R u_h \times \nu_F$ on boundary facets $F \in \mathcal{F}(\partial \Omega)$.
This and Lemma \ref{lem:local-data-oscillation} establish the existence of some mesh-independent constant $\Ceff>0$ with $\Ceff^{-1}\etares(\T)\leq \trb{u-Ru_h}_\pw + \osc q(f,\T)$ for arbitrary $q\in\mathbb N_0$. This concludes the proof.\qed
\end{proof}

While the focus of this paper is on the HHO methodology, the framework of Section \ref{sec:blocks} also applies to other skeletal methods as well.	
The following example covers a hybridized discontinuous Galerkin (HDG) FEM from \cite{Oikawa2015} with the Lehrenfeld-Sch\"oberl stabilization.
\begin{example}[HDG FEM]
	Let $V_h \coloneqq P_{k+1}(\mathcal{T}) \times P_k(\mathcal{F}(\Omega))$ for $k \in \mathbb{N}_0$. An equivalent formulation to the HDG FEM from \cite{Oikawa2015} seeks $u_h \in V_h$ with
	\begin{align*}
		a_\pw(R u_h, R v_h) + s_h(u_h,v_h) = (f, v_\mathcal{T}) \quad\text{for any } v_h = (v_\mathcal{T},v_\mathcal{F}) \in V_h.
	\end{align*}
	Here, $R : V_h \to P_{k+1}(\mathcal{T})$ is defined as in \eqref{eqn:R} and
	\begin{align*}
		s_h(v_h,w_h) \coloneqq \sum_{T\in\T}\sum_{F\in\F(T)} h_F^{-1}\langle \Pi_{F,k}(v_\mathcal{T}|_T - v_\mathcal{F}|_F), w_\mathcal{T}|_T - w_\mathcal{F}|_F \rangle_{L^2(F)}
	\end{align*}
	for any $v_h = (v_\mathcal{T}, v_\mathcal{F}), w_h = (w_\mathcal{T}, w_\mathcal{F}) \in V_h$.
	This method is also known under the label of weak Galerkin FEM \cite{WangYe2013}.
	It is straightforward to verify that $\nabla_\pw R u_h$ satisfies \eqref{eq:solution_property}--\eqref{eq:div_0_property}. Notice that \eqref{eq:solution_property} also holds for $k = 0$ without any modification.
	Therefore, Theorem \ref{thm:reliability} leads to the reliable a~posteriori estimate
	\begin{align*}
		&\trb{u-R u_h}_\pw^2 \leq \eta_{\mathrm{res}}^2(\mathcal{T}) \coloneqq \|h_\T(f + \Delta_{\rm pw} Ru_h)\|^2\\
		&~ + \sum_{F\in \cF(\Omega)} \ell(F)\|[\nablapw Ru_h ]_F\cdot \nu_F \|_{L^2(F)}^2 + \sum_{F\in \cF} \ell(F)\|[\nablapw Ru_h ]_F\times \nu_F \|_{L^2(F)}^2.
	\end{align*}
	The efficiency $\eta_{\mathrm{res}}(\mathcal{T}) \lesssim \trb{u-R u_h}_\pw$ follows from the arguments in the proofs of Lemma \ref{lem:efficiency}--\ref{lem:local-data-oscillation}.
\end{example}
\section{Equilibrium-based a posteriori HHO error analysis}
\label{sec:eq}

\newcommand{\gRuh}{\nabla_{\rm pw} Ru_h}
\newcommand{\qdelta}{q^\Delta}
\newcommand{\tf}{\tilde f}

\newcommand{\supp}{\text{supp}} 
\newcommand{\mz}{|\partial \oz \cap \partial \Omega|} 
\newcommand{\RTz}{
RT_{k+1}^0(\mathcal T (z)) %
}
\newcommand{\RTzk}{
RT_{k}^0(\mathcal T (z)) %
}
\newcommand{\RTzkp}{
RT_{k+1}^0(\mathcal T (z)) %
}
\newcommand{\PiRT}{\Pi_{\RTzk}}
\renewcommand{\PiRT}{\Pi_{RT}^k}
\newcommand{\PiRTk}{\Pi_{\RTzk}}
\newcommand{\Tz}{\mathcal T (z)}
\newcommand{\Lz}{L^2(\oz)}
\newcommand{\niLz}{_{\Lz}}
\newcommand{\niLT}{_{L^2(T)}}
\newcommand{\nLz}[1]{||#1||\niLz}
\newcommand{\RTkppwz}{RT^{k+1}_{{\rm pw},0}(\mathcal T (z))}

The residual-based guaranteed upper bound (GUB) of the error $\trb{u - R u_h}_\pw$ from Subsection \ref{sub:HHO_eta} employs explicit constants that may lead to overestimation in higher dimensions and for different triangular shapes.
This section utilizes equilibrated flux reconstructions \cite{Ainsworth2005,Ain:07,ErnVohralik2015,bertrand_weakly_2019,bertrand_opt} to establish, up to the well-known Poincar\'e constant $C_P \leq 1/\pi$, a constant-free guaranteed upper bounds for a tight error control.

\subsection{Guaranteed error control}
\label{sub:Motivation_eq}
The guaranteed upper bounds of this section involves \emph{two post-processings} of the potential reconstruction $R u_h \in P_{k+1}(\T)$ of the discrete solution $u_h$ to \eqref{eqn:HHO}.
First, the patch-wise design of a flux reconstruction $Q_p \in RT_{k+p}(\T)$ with $p \in \mathbb{N}_0$ from \cite{AinOde:93,BraPilSch:09,ErnVohralik2020} provides an $H(\mathrm{div},\Omega)$-conforming approximation to $\nablapw R u_h$ with the equilibrium $\Pi_{r} f + \div Q_p = 0$ in $\Omega$
and $r$ from \eqref{def:r} below.
Second, the nodal average $\mathcal{A} R u_h \in S^{k+1}_0(\T)\subset V$ results in an $V$-conforming approximation of $R u_h$ by averaging all values of the discontinuous function $R u_h$ at each Lagrange point of $S^{k+1}_0(\T)$. %
This, the split \eqref{eqn:decomp}, and the solution property \eqref{eqn:HHO_solution_property} give rise to the guaranteed upper bound (GUB)
\begin{align}
\label{eq:definitionetaq}
\eta_{\mathrm{eq},p}^2(\T) := \left(C_P\osc{r}(f,\T) + \| \QQ_p- \nablapw R u_h\|\right)^2 + \trb{(1 - \mathcal{A})R u_h}_\pw^2
\end{align}
with $r \in \mathbb{N}_0$ defined by
\begin{align}\label{def:r}
    r \coloneqq 0 \text{ if } k = 0 \quad\text{and}\quad r \coloneqq k + p \text{ if } k \geq 1.
\end{align}
The main result of this section states the reliability and efficiency (up to data oscillations) of $\eta_{\mathrm{eq},p}$ for all parameters $p \in \mathbb{N}_0$.
\begin{theorem}[equilibrium-based GUB for HHO]\label{thm:GUB-HHO}
    Let $u \in V$ resp.\ $u_h \in V_h$ solve \eqref{eq:poisson} resp.\ \eqref{eqn:HHO}.
    Given a parameter $p \in \mathbb{N}_0$,
    there exists $\QQ_p\in RT_{k+p}(\mathcal{T})$ such that
    the error estimator $\eta_{\mathrm{eq},p}(\T)$ from \eqref{eq:definitionetaq}
    is an efficient GUB 
    \begin{align}\label{ineq:equilibrium-equivalence}
        \trb{u-R u_h}_\pw \le\eta_{\mathrm{eq},p}(\T) \leq \cnst{cnst:equilibrium-efficiency}\big(\trb{u-R u_h}_\pw + \osc{q}(f, \mathcal{T})\big).
    \end{align}
    for any $q\in \mathbb N_0$ and $\osc{k-1}(f, \mathcal{T}) \leq\cnst{cnst:equilibrium-reliability}\eta_{\mathrm{eq},p}(\T)$.
    The constants $\newcnst\label{cnst:equilibrium-efficiency}$ and $\newcnst\label{cnst:equilibrium-reliability}$ exclusively depend on the polynomial degree $k \in \mathbb{N}_0$, the parameter $q \in \mathbb{N}_0$, and the shape-regularity of $\T$.
\end{theorem}
At least two technical contributions for the proof of Theorem~\ref{thm:GUB-HHO} are of broader interest.
A first contribution to the HHO literature
is the local equivalence
of the original HHO stabilization $s_h$ from \eqref{eqn:s_h0} and the alternative stabilization $\tilde{s}_h(v_h,v_h) \coloneqq \sum_{T \in \T} \tilde{s}_T(v_h,v_h)$ from \cite{DiPietroDroniou2020} defined, for  $v_h=(v_\T, v_\F)\in V_h$, by
\begin{align}
    \label{def:s-tilde}
    \tilde{s}_T(v_h,v_h)\coloneqq h_T^{-2} \|\Pi_{T,k} (v_\T &- R v_h)\|_{L^2(T)}^2 \\&+ \sum_{F \in \F(T)} h_F^{-1} \|\Pi_{F,k} (v_\F - R v_{h}|_{T})\|_{L^2(F)}^2.\nonumber
\end{align}
A second result of separate interest in the HHO literature (cf.~\cite{DiPietroDroniou2020} where the efficiency in \eqref{eqn:efficient_stabilization} is left open) is the efficiency of the stabilizations from Theorems \ref{thm:equivalence-stabilization}--\ref{thm:best-approximation} below,
\begin{align}\label{eqn:efficient_stabilization}
    \tilde{s}_h(v_h,v_h)^{1/2}\approx s_h(u_h,u_h)^{1/2} \lesssim \trb{u - R u_h}_\pw + \osc{k+p}(f,\part).
\end{align}

The subsequent subsection continues with some explanations on the flux reconstruction $\QQ_p\in RT_{k+p}(\mathcal{T})$ that is defined by local minimization problems on each vertex patch.
Appendix A complements the discussion with an algorithmic two-step procedure for the computation of $\QQ_{p}-\nablapw R u_h$ in 2D.
The efficiency of the averaging $\trb{(1 - \mathcal{A})R u_h}_\pw$ up to data oscillations follows in Subsections \ref{sub:local-equivalence-of-stabilizations}--\ref{sub:efficiency_stabilization}.
Subsection \ref{sub:proof-efficiency-GUB} concludes with  the proof of Theorem \ref{thm:GUB-HHO}.

\subsection{Construction of equilibrated flux}
\renewcommand{\PiRT}{\Pi_{RT_k^{\rm pw}}}
\label{sub:Computation_Q}
This subsection defines the post-processed $H(\mathrm{div},\Omega)$-conforming equilibrated flux $\QQ_p\in RT_{k+p}(\mathcal{T})$ that enters the GUB $\eta_{\mathrm{eq},p}$ from \eqref{eq:definitionetaq} based on local patch-wise minimization problems in the spirit of \cite{BraPilSch:09,ErnVohralik2015,ErnVohralik2020}.

Consider the shape-regular vertex-patch $\omega(z)\coloneqq \mathrm{int}(\bigcup\{T\in\T(z)\})$ covered by the neighbouring simplices $\T(z)\coloneqq\{T\in\T : z\in T\}$ sharing a given vertex $z \in \mathcal{V}$ with the facet spider $\F(z)\coloneqq\{F\in\F\ :\ z\in F\}$.
Recall the space of piecewise Raviart-Thomas functions $RT_{k}^{\pw}(\T)$ from Subsection \ref{sub:notation} and define
\begin{align*}
    L^2_0(\omega(z)) &\coloneqq \{f\in L^2(\omega(z)) : (f,1)\niLz = 0\},\\
    L^2_*(\omega(z)) &\coloneqq \begin{cases}
    L^2_0(\omega(z))& \text{ if } z \in \mathcal{V}(\Omega),\\
    L^2(\oz) &\text{ else},
    \end{cases}\\
    H^1_*(\omega(z)) &\coloneqq 
    \begin{cases}
        H^1(\omega(z)) \cap L^2_0(\omega(z))\quad\text{ if } z \in \mathcal{V}(\Omega),\\
        \{v \in H^1(\omega(z)) : v = 0 \text{ on } \partial \Omega \cap \bigcup\F(z)\}\quad\text{ else},
    \end{cases}\\
    \Hdivz &\coloneqq
    \begin{cases}
        \{r \in H(\mathrm{div},\omega(z)) : r \cdot \nu = 0 \text{ on } \partial \omega(z)\} \quad\text{ if } z \in \mathcal{V}(\Omega),\\
        \{r \in H(\mathrm{div},\omega(z)) : r \cdot \nu = 0 \text{ on } \partial \omega(z) \setminus \bigcup\F(z)\}\quad\text{ else},
    \end{cases}\\
    RT_k^0(\T(z))&\coloneqq\;RT_{k}^{\pw}(\T(z)) \cap \Hdivz.
\end{align*}
Throughout the remaining parts of this section, abbreviate $\GG \coloneqq \nabla_\pw R u_h \in P_k(\mathcal{T};\R^n)$.
Given a vertex $z\in\mathcal V$ with the $P_1$-conforming nodal basis function $\varphi_z\in S^1(\T)$,
the property \eqref{eqn:HHO_solution_property} provides compatible data
\begin{align}\label{def:f-z}
    f_z \coloneqq\,
    \begin{cases}
        \Pi_{p}(\varphi_z \Pi_0 f- \GG \cdot \nabla \varphi_z) &\mbox{if } k = 0,\\
        \Pi_{k+p}(\varphi_z f- \GG \cdot \nabla \varphi_z) &\mbox{if } k \geq 1
    \end{cases}\quad\in L^2_*(\oz)
\end{align}
such that the discrete affine space \begin{align}\label{eqn:Qhz_def}
    \mathcal{Q}_h(z) \coloneqq \{\tau_z \in RT^0_{k+p}(\T(z))\ :\ \mathrm{div}\, \tau_z + f_z = 0 \text{ in } \Omega\} \neq \emptyset
\end{align}
is not empty. Consequently,
\begin{align}
\label{eq:discrete_minimisation}
\Rzh \coloneqq \argmin\limits_{\tau_z \in \mathcal{Q}_h(z)}
     \| \tau_z - \mathcal{I}_{RT} (\varphi_z \GG)\|\niLz = \Pi_{\mathcal{Q}_h(z)} \mathcal{I}_{RT} (\varphi_z \GG)
\end{align}
is well defined as the $L^2$ projection $\Pi_{\mathcal{Q}_h(z)} \mathcal{I}_{RT} (\varphi_z \GG)$ of $\mathcal{I}_{RT} (\varphi_z \GG)$ onto $\mathcal{Q}_h(z)$
with the piecewise Raviart-Thomas interpolation $\mathcal{I}_{RT} : H^1(\T;\R^n)\to RT^\pw_{k+p}(\T)$ \cite[Section III.3.1]{boffi_mixed_2013}.
In the case $p \geq 1$, $\varphi_z \GG \in P_{k+1}(\mathcal{T}(z)) \subset RT_{k+p}^\pw(\mathcal{T}(z))$
is a piecewise Raviart-Thomas function of degree $k+p$.
Hence $\mathcal{I}_{RT}(\varphi_z \GG) = \varphi_z \GG$ and $\mathcal{I}_{RT}$
could be omitted in the formula \eqref{eq:discrete_minimisation}.
The partition of unity $\sum_{z\in\mathcal V}\varphi_z\equiv 1$ and $\GG = \mathcal{I}_{RT} \GG = \sum_{z \in \mathcal{V}} \mathcal{I}_{RT} (\varphi_z \GG)$ show that the sum $\QQ_p = \sum_{z \in \mathcal V} \Rzh\in H(\ddiv, \Omega)$ of the patch-wise contributions 
satisfies
\begin{align}
    &\ddiv \QQ_p = 
    \begin{cases}
        \sum\limits_{z \in \mathcal V} \Pi_{p} (\GG \cdot\nabla \varphi_z - \varphi_z \Pi_0 f) = -\Pi_{0} f &\mbox{if } k = 0,\\
        \sum\limits_{z \in \mathcal V} \Pi_{k+p} (\GG \cdot\nabla \varphi_z - \varphi_z f) = -\Pi_{k+p} f &\mbox{if } k \geq 1,
    \end{cases}\label{eqn:div_Q_f}\\
    &\|\QQ_{p}-\GG \|_{L^2(\Omega)}^2 \lesssim \sum_{z\in\mathcal V}
    \|\Rzh - %
    \mathcal{I}_{RT}(\varphi_z \GG)
    \|_{L^2(\oz)}^2.\label{eqn:Q_q_small}
\end{align}
This establishes the flux reconstruction $Q_p$.
The efficiency of the flux reconstruction will be based on the following general equivalence.
\newcommand{\Cs}{C_{\rm s}}
\begin{lemma}[control of $H(\ddiv)$ minimization by residual {\cite{BraPilSch:09,ErnVohralik2020}}]\label{lem:uniform_stability_mixed}
\newcommand{\sz}{\sigma_z}
Given any vertex $z\in\mathcal V$, a piecewise Raviart-Thomas function $\sz \in RT_{q}^{{\rm pw}}(\T(z))$ and a piecewise polynomial $r_z \in {P}_{q}\left(\mathcal{T}(z)\right)$ of degree $q \in \mathbb{N}_0$, define the residual
\begin{align}\label{eq:compatibility}
 \mathrm{Res}_z(v)\coloneqq \sum_{T \in \mathcal{T}(z)}\Big(
\left(r_{z}, v\right)_{L^2(T)}+\left\langle\sz\cdot \nu_T, v\right\rangle_{L^2(\partial T)}\Big)
\end{align}
for all $v\in H^1(\Omega)$. If $z\in \mathcal V(\Omega)$ is an interior vertex, then suppose additionally that
$\mathrm{Res}_z(1) = 0$.
Then
\begin{align}
\label{lem:uniform_stability_mixedineq}
\min\limits_{\substack{\tau_z \in RT^0_q(\mathcal{T}(z))\\\mathrm{div} \tau_z = r_z + \mathrm{div}_\pw \sigma_z}}
\left\|\tau_z-\sz\right\|_{L^2(\oz)} 
\leq \Cs
\max\limits_{%
       \substack{%
         v \in H^{1}_*(\omega(z))
         \\
         \|\nabla v\|_{L^2(\oz)}=1
       }
     } \mathrm{Res}_z(v)
\end{align}
holds for a constant $\Cs$ that exclusively depends on the shape-regularity (and is in particular independent of the polynomial degree $q$).
\end{lemma}
\begin{proof}
    The assertion follows from \cite[Theorem 7]{BraPilSch:09} in $n=2$ dimensions and \cite[Corollaries 3.3, 3.6, and 3.8]{ErnVohralik2020} in $n=3$ dimensions.\qed
\end{proof}

\begin{remark}\label{rem:local_osc_eff}
The patch-wise construction of $Q_p$ in Subsection \ref{sub:Computation_Q} typically 
generates local data oscillation $\osc{k+p}(\varphi_z f,\mathcal{T}(z))$ in the error analysis as in the proof of Theorem \ref{thm:GUB-HHO} in Subsection \ref{sub:proof-efficiency-GUB} below or, e.g.,  \cite[Theorem 3.17]{ErnVohralik2015}.
A straightforward computation $\osc{k+p}(\varphi_z f,\mathcal{T}(z)) \leq \osc{k+p-1}(f,\mathcal{T}(z))$ apparently leads to a loss of one degree
in the data oscillation but Lemma \ref{lem:local-data-oscillation} verifies 
\begin{align}\label{eqn:local_osc_eff}
    \osc{k+p}^2(\varphi_z f, \part(z))\lesssim \min_{v_{h}\in P_{k+1}(\T(z))} \|\nablapw(u - v_{h})\|_{L^2(\omega(z)}^2 + \osc{q}^2(f,\part(z))
\end{align}
for any $p, q\in\N_0$.
This allows for efficiency of the data oscillations 
on the right-hand side of the efficiency estimate \cite[Formula (3.42)]{ErnVohralik2015} and leads to a corresponding refinement in 
\cite[Theorem 3.17]{ErnVohralik2015}.
\end{remark}

\subsection{Local equivalence of stabilizations}
\label{sub:local-equivalence-of-stabilizations}
The first improvement to the current HHO literature is the local equivalence of the two stabilizations $\tilde{s}_h$ from \eqref{def:s-tilde} and $s_h$ from \eqref{eqn:s_h0}.
The authors of this paper could not find any motivation for the alternative stabilization $\tilde{s}_h$ in the error analysis of \cite[Section 4]{DiPietroDroniou2020} and suggest to apply Theorem \ref{thm:equivalence-stabilization} below to \cite[Theorem 4.7]{DiPietroDroniou2020} to recover the results therein for the original HHO stabilization $s_h$.
Recall the local stabilization $\tilde{s}_T$ in $\tilde{s}_h(v_h,v_h) \coloneqq \sum_{T \in \T} \tilde{s}_T(v_h,v_h)$ from \eqref{def:s-tilde}
and $\SpTF v_h =
\Pi_{F,k}
\left(
v_{\T} +
(1-\Pi_{T,k})R  v_{h} 
\right)|_{T}-v_{\F}|_{F}$ in the definition of $s_h$ from \eqref{eqn:s_h0}.
\begin{theorem}[local equivalence of stabilizations]\label{thm:equivalence-stabilization}
	Any $v_h = (v_{\part},v_\mathcal{F}) \in V_h$ and $T \in \part$ satisfy
	\begin{align}
		\cnst{cnst:equivalence-norms-LHS}^{-1}\tilde{s}_T(v_h,v_h) \leq \sum_{F \in \mathcal{F}(T)} h_F^{-1} \|\SpTF v_h\|_{L^2(F)}^2 \leq \cnst{cnst:equivalence-norms-RHS}\tilde{s}_T(v_h,v_h).
		\label{ineq:efficiency-trace}
	\end{align}
	The constants $\newcnst\label{cnst:equivalence-norms-LHS}$ and $\newcnst\label{cnst:equivalence-norms-RHS}$ exclusively depend on the polynomial degree $k$ and the shape regularity of $\mathcal{T}$.
\end{theorem}
\begin{proof}
    The second inequality in \eqref{ineq:efficiency-trace} follows directly from a triangle inequality and an inverse estimate. Therefore, the proof focuses on the first inequality in \eqref{ineq:efficiency-trace}.
	Given $v_h = (v_\part,v_\mathcal{F}) \in V_h$ and $T \in \part$, let $\varphi_k \coloneqq (\Pi_k R v_h - v_\T)|_{T} \in P_k(T)$.
	Since $\SpTF v_h = \Pi_{F,k}(R v_{h|T} - v_\F|_F - \varphi_k)$, the triangle inequality $\|\Pi_{F,k}(R v_{h|T} - v_{\F})\|_{L^2(F)} \leq \|\SpTF v_h\|_{L^2(F)} + \|\varphi_k\|_{L^2(F)}$, the discrete trace inequality $\|\varphi_k\|_{L^2(F)} \lesssim h_F^{-1/2}\|\varphi_k\|_{L^2(T)}$,
	and the shape-regularity $h_F \approx h_T$ for all $F \in \mathcal{F}(T)$ reveal
	\begin{align}
		&\sum_{F \in \mathcal{F}(T)} h_F^{-1}\|\Pi_{F,k}(R v_{h|T} - v_{\F})\|_{L^2(F)}^2\nonumber\\
		&\qquad\qquad\lesssim \sum_{F \in \mathcal{F}(T)} h_F^{-1} \|\SpTF v_h\|_{L^2(F)}^2 + h_T^{-2}\|\varphi_k\|_{L^2(T)}^2.
		\label{ineq:proof-efficiency-trace-triangle-inequality}
	\end{align}
	Since $\Pi_0 \varphi_k = 0$ (from the design of $Rv_h$), a Poincar\'e inequality shows
	\begin{align}
		h_T^{-2}\|\varphi_k\|_{L^2(T)}^2 \leq C_P^2\|\nabla \varphi_k\|_{L^2(T)}^2.
		\label{ineq:proof-efficiency-trace-Poincare}
	\end{align}
	On the one hand, an integration by parts provides
	\begin{align}
		\|\nabla \varphi_k\|_{L^2(T)}^2 = - (\Pi_k R v_h - v_\T, \Delta \varphi_k)_{L^2(T)} + \langle\varphi_k, \nabla \varphi_k \cdot \nu_T\rangle_{L^2(\partial T)}.
		\label{eq:proof-upper-bound-trace-ibp-1}
	\end{align}
	On the other hand, an integration by parts and the definition of $R$ imply
	\begin{align}
		&- (R v_h, \Delta \varphi_k)_{L^2(T)} = (\nabla R v_h, \nabla \varphi_k)_{L^2(T)} - \langle R v_{h|T}, \nabla \varphi_k \cdot \nu_T\rangle_{L^2(\partial T)}\nonumber\\
		&\qquad= - (v_{\T}, \Delta \varphi_k)_{L^2(T)} + \sum_{F \in \mathcal{F}(T)} \langle v_{\F} - R v_{h|T}, \nabla \varphi_k \cdot \nu_T\rangle_{L^2(F)}.
		\label{eq:proof-upper-bound-trace-ibp-2}
	\end{align}
	Since $\Delta \varphi_k \in P_k(T)$, the $L^2$ projection $\Pi_k$ on the right-hand side of \eqref{eq:proof-upper-bound-trace-ibp-1} is redundant. Hence, the combination of \eqref{eq:proof-upper-bound-trace-ibp-1}--\eqref{eq:proof-upper-bound-trace-ibp-2} with $\nabla \varphi_k \cdot \nu_{T|F} \in P_k(F)$ for all $F \in \mathcal{F}(T)$
	results in
	\begin{align}
        \|\nabla \varphi_k\|_{L^2(T)}^2 = \sum_{F \in \mathcal{F}(T)} \langle \Pi_{F,k}(v_{\F} - R v_{h|T} + \varphi_k), \nabla \varphi_k \cdot \nu_T \rangle_{L^2(F)}.
        \label{eq:proof-upper-bound-trace-ibp-3}
	\end{align}
	A Cauchy inequality on the right-hand side of \eqref{eq:proof-upper-bound-trace-ibp-3}, a discrete trace inequality, and $\SpTF v_h = \Pi_{F,k}(R v_{h|T} - v_{\F|F} - \varphi_k)$ for all $F \in \mathcal{F}(T)$ lead to
	\begin{align}
		\|\nabla \varphi_k\|_{L^2(T)}^2 \lesssim \sum_{F \in \mathcal{F}(T)} h_F^{-1} \|\SpTF v_h\|_{L^2(F)}^2.
		\label{ineq:bound-gradient-phi-k}
	\end{align}
	Since $\tilde{s}_T(v_h,v_h) = \sum_{F \in \mathcal{F}(T)} h_F^{-1}\|\Pi_{F,k}(R v_{h|T} - v_{\F})\|_{L^2(F)}^2 + h_T^{-2}\|\varphi_k\|_{L^2(T)}^2$,
	the combination of \eqref{ineq:proof-efficiency-trace-triangle-inequality}--\eqref{ineq:proof-efficiency-trace-Poincare} with \eqref{ineq:bound-gradient-phi-k} concludes the proof of \eqref{ineq:efficiency-trace}.\qed
\end{proof}

\subsection{Efficiency of the stabilization}\label{sub:efficiency_stabilization}
The second improvement to the HHO literature
is a quasi-best approximation estimate along the lines of the seminal paper \cite[Theorem 4.10]{ErnZanotti2020}.
In combination with Theorem \ref{thm:equivalence-stabilization}, this, in particular, provides the efficiency \eqref{eqn:efficient_stabilization} of the stabilization up to data oscillation.

\begin{theorem}[quasi-best approximation up to data oscillation]\label{thm:best-approximation}
    For any $p\in\mathbb N_0$, the solution $u$ to \eqref{eq:poisson} and the discrete solution $u_h$ to  \eqref{eqn:HHO} satisfy
	\begin{align*}
		\trb{u - R u_h}_\pw &+ s_h(u_h,u_h)^{1/2}\\
		&\leq \cnst{cnst:best-approximation}\Big(\min_{v_{k+1} \in P_{k+1}(\part)} \trb{u -  v_{k+1}}_\pw + \osc{k+p}(f,\part)\Big).
	\end{align*}
	The constant $\newcnst\label{cnst:best-approximation}$ exclusively depends on $k$, $p$, and the shape regularity of $\mathcal{T}$.
\end{theorem}
\begin{proof}
    Given $k,p \in \mathbb{N}_0$, let $\widetilde{u} \in V$ solve the Poisson model problem $-\Delta \widetilde{u} = \Pi_{k+p} f$ with the right-hand side $\Pi_{k+p} f$.
    Subsection 4.3 in \cite{ErnZanotti2020} constructs a stable enriching operator $\mathcal{J} : V_h \to V$ with local bubble-functions from \cite{Ver:13}.
    A modification, where the polynomial degree $k-1$ in \cite[Eq.~(4.16)]{ErnZanotti2020} is replaced by $k+p$, leads to a right-inverse $\mathcal{J} : V_h \to V$ of the interpolation $\mathrm{I}: V_h \to V$ with the stability $\trb{\mathcal{J} v_h}^2 \lesssim a_h(v_h,v_h)$ and the additional $L^2$ orthogonality $\mathcal{J} v_h - v_\mathcal{T} \perp P_{k+p}(\mathcal{T})$ for all $v_h = (v_\mathcal{T},v_\mathcal{F}) \in V_h$.
    The extra orthogonality allows for
$$
        (\Pi_{k+p} f, \mathcal{J} v_h)_{L^2(\Omega)} = (\Pi_{k+p} f, v_\mathcal{T})_{L^2(\Omega)} = (f,v_\mathcal{T})_{L^2(\Omega)}\quad\text{for all } (v_\T, v_\F)\in V_h.
    $$
    Consequently, the smoother $\mathcal{J}$ leads to a discrete solution $u_h=(u_\T,u_\F)\in V_h$ in the modified HHO discretization of \cite{ErnZanotti2020} as a quasi-best approximation of the above solution $\tilde u$. The point is that $u_h\in V_h$ coincides with the original HHO solution $u_h$ from \eqref{eqn:HHO}.
    The arguments from the proof of \cite[Theorem 4.10]{ErnZanotti2020} reveal the quasi-best approximation 
    $$\trb{\tilde{u} - R u_h}_\pw + s_h(u_h,u_h)^{1/2} \lesssim \min_{v_{k+1} \in P_{k+1}(\T)} \trb{\tilde{u} - v_{k+1}}_\pw$$ also for the above modified smoother $\mathcal{J}$.
    This, the triangle inequalities $\trb{u - R u_h}_\pw \leq \trb{u - \tilde{u}} + \trb{\tilde{u} - R u_h}_\pw$ and
    $\trb{\tilde{u} - v_{k+1}}_\pw \leq \trb{u - \tilde{u}} + \trb{u - v_{k+1}}_\pw$ for any $v_{k+1} \in P_{k+1}(\mathcal{T})$, and the standard perturbation bound $\trb{u - \tilde{u}} \leq C_P\osc{k+p}(f,\T)$ conclude the proof. \qed
\end{proof}

\subsection{Stabilization-free efficiency of averaging}\label{sub:averaging-efficiency}
The main result of this subsection establishes the stabilization-free efficiency of the nodal averaging technique.
\begin{theorem}[averaging is efficient]\label{thm:averaging-efficiency}
    Let $u \in V$ resp.\ $u_h \in V_h$ solve \eqref{eq:poisson} resp.\  \eqref{eqn:HHO}. Then $R u_h$ and $\mathcal{A} R u_h$ satisfy, for any $p \in \mathbb N_0$, that
    \begin{align}
        \cnst{cnst:efficiency-averaging}^{-1}\trb{(1 - \mathcal{A}) R u_h}_\pw 
        &\leq \trb{u - R u_h}_\pw + \osc{k+p}(f,\part).
        \label{ineq:averaging-efficiency-oscillation}
    \end{align}
    The constant $\newcnst\label{cnst:efficiency-averaging}$ exclusively depends on $k$, $p$, and the shape regularity of the triangulation $\mathcal{T}$.
\end{theorem}
The proof can follow the proof of \cite[Theorem 4.7]{DiPietroDroniou2020} but additionally utilizes the two significant improvements from Subsections
\ref{sub:local-equivalence-of-stabilizations}--\ref{sub:efficiency_stabilization} that allow a stabilization-free efficiency in \eqref{ineq:averaging-efficiency-oscillation}.
\begin{proof}
	Theorem 4.7 in \cite{DiPietroDroniou2020} shows $\trb{(1 - \mathcal{A})R u_h}_\pw^2 \lesssim \trb{u - R u_h}_\pw^2 + \tilde{s}_h(u_h,u_h)$. This, the equivalence $\tilde{s}_h(u_h,u_h) \approx s_h(u_h,u_h)$ of stabilizations (from Theorem \ref{thm:equivalence-stabilization}), and 
	the efficiency $s_h(u_h,u_h) \lesssim \trb{u - R u_h}_\pw^2 + \osc{k+p}^2(f,\mathcal{T})$ (from Theorem \ref{thm:best-approximation}) imply the assertion. \qed
\end{proof}

\begin{remark}[$p$-robustness]
    The $H^1(\Omega)$-conforming post-processing of $R u_h$ from \cite{ErnVohralik2020} provides an efficiency constant independent of the polynomial degree $k$.
    The right-hand side of \cite[Corollary 4.2]{ErnVohralik2020} involves the stabilization-related term $\sum_{F \in \F}h_F^{-1}\|\Pi_{F,0}[R u_h]_F\|_{L^2(F)}^2$.
    It remains an open question whether this term can be controlled $p$-robustly by $\trb{u - R u_h}_\pw + \osc{k}(f,\T)$ (with a multiplicative 
    constant independent of the polynomial degree $k$).
\end{remark}

\subsection{Proof of Theorem \ref{thm:GUB-HHO}}\label{sub:proof-efficiency-GUB}
Let $p \in \mathbb{N}_0$ and $r = 0$ if $k = 0$ and $r = k+p$ if $k \geq 1$ as in \eqref{def:r} be given. Recall the abbreviation $\GG = \nabla_\pw R u_h \in P_k(\mathcal{T})$ with the discrete solution $u_h \in V_h$ to \eqref{eqn:HHO} and let $\QQ_p = \sum_{z \in \mathcal V} \Rzh\in H(\ddiv, \Omega)$ denote the flux reconstruction from Subsection \ref{sub:Computation_Q}. The proof establishes \eqref{ineq:equilibrium-equivalence} in five steps.\medskip

\emph{Step 1} provides the GUB $\trb{u - R u_h}_\pw \leq \eta_{\mathrm{eq},p}(\mathcal{T})$.
This can follow from the paradigm of \cite{Ainsworth2005,Ain:07,ErnVohralik2015} as outlined below.
The choice $\QtoG \coloneqq \GG$ and $w \coloneqq \mathcal{A} R u_h$ in \eqref{eqn:decomposition_full} and a triangle inequality lead to
\begin{align*}
    \trb{u - R u_h}_\pw^2
    \leq (\trb{f + \ddiv Q_p}_* + \trb{\mathrm{div}(Q_p - \GG)}_*)^2 + \trb{(1 - \mathcal{A}) R u_h}_\pw^2.%
\end{align*}
Since $\ddiv Q_p + \Pi_{r} f = 0$ vanishes in $\Omega$ by \eqref{eqn:div_Q_f}, a piecewise Poincar\'e inequality shows $\trb{f + \ddiv Q_p}_* \leq C_P\osc{r}(f,\T)$. This, the bound $\trb{\mathrm{div}(Q_p - \GG)}_* \leq \|Q_p - \GG\|$ from the definition of $\trb{\bullet}_*$, and the previously displayed formula result in the reliability $\trb{u - R u_h}_\pw \leq \eta_{\mathrm{eq,p}}(\mathcal{T})$.\medskip

\emph{Step 2} establishes $\osc{k-1}(f,\mathcal{T}) \lesssim \eta_{\mathrm{eq},p}(\mathcal{T})$. 
Lemma \ref{lem:local-data-oscillation} provides 
\begin{align}
    \osc{k-1}(f,\mathcal{T})\lesssim \trb{u - R u_h}_\pw +\osc r(f, \T).
\end{align}
This, Step 1 and $\osc r(f, \T)\lesssim\etaeqnew(\T)$ from \eqref{eq:definitionetaq} conclude the proof of Step 2.\medskip

\newcommand{\IRT}{\mathcal I_{\RT}}
\emph{Step 3} reveals, for any $q\in\N_0$, the efficiency of the flux reconstruction
\begin{align}\label{ineq:proof-GUB-efficiency-flux}
    \|\QQ_p - \GG\| \lesssim \trb{u-R u_h}_\pw + \osc{q}(\part, f)
\end{align}
for any polynomial degree $k \geq 1$.
The case $k=0$ follows in Step 4 below.
Recall $f_z = \Pi_{k+p}(\varphi_z f - \GG\cdot \nabla \varphi_z)$ from \eqref{def:f-z} for any vertex $z\in\mathcal V$ in the construction of $Q_p$ from Subsection \ref{sub:Computation_Q} and set
$\sigma_z \coloneqq \IRT(\varphi_z \GG)\in RT_{k+p}^{{\rm pw}}(\T(z))$. %
The piecewise Raviart-Thomas interpolation $\IRT:H^1(\T;\R^n)\to RT_{k+p}^{\rm pw}(\T)$ satisfies the well-known commuting diagram properties \cite[Section 2.5.1]{boffi_mixed_2013}
\begin{align}\label{eq:commuting-diagram-I-RT}
    \divpw\circ\IRT = \Pi_{k+p}\circ\divpw&&\text{and}&&\gamma_T\circ\IRT|_F=\Pi_{F,k+p}\circ \gamma_T
\end{align}
for any facet $F\in\F(T)$ of a simplex $T\in\T$ and the normal trace $\gamma_T:H^1(T;\R^n)\to L^2(\partial T)$ with $\gamma_T \sigma \coloneqq \sigma \cdot \nu_T$ for $\sigma\in H^1(T;\R^n)$.
This and elementary algebra with the product rule $\divpw (\varphi_z \GG)=\varphi_z\divpw \GG + \nabla \varphi_z\cdot \GG \in P_k(\mathcal{T}(z))$ imply
\begin{align}\label{ineq:proof-GUB-efficiency-p=0-divergence}
    r_z \coloneqq -\mathrm{div}_\pw \sigma_z - f_z = - \varphi_z \mathrm{div}_\pw \GG-  \Pi_{k+p}(\varphi_zf)\in P_{k+p}(\T(z)).
\end{align}
Recall the residual $\mathrm{Res}_z(v)$ with $v\in H^1(\Omega)$ for the vertex $z\in\mathcal V$ from \eqref{eq:compatibility}. The commuting diagram property \eqref{eq:commuting-diagram-I-RT}--\eqref{ineq:proof-GUB-efficiency-p=0-divergence} establish the identity 
\begin{align*}
    \mathrm{Res}_z(1) = -(\mathrm{div}_\pw \GG + f, \varphi_z)_{L^2(\omega(z))} + \sum_{T \in \mathcal{T}(z)} \langle \GG \cdot \nu_T, \varphi_z\rangle_{L^2(\partial T)}.
\end{align*}
This, a piecewise integration by parts, and the property \eqref{eqn:HHO_solution_property} verify $\mathrm{Res}(1) = (\GG,\nabla \varphi_z)_{L^2(\omega(z))} - (f,\varphi_z)_{L^2(\omega(z))}=0$ for any interior vertex $z \in \mathcal{V}(\Omega)$.
Hence, Lemma \ref{lem:uniform_stability_mixed} applies for any vertex $z\in\mathcal V$ and provides
\begin{align}\label{eq:efficiencyHdiv}
    \|\Rzh - \mathcal{I}_{RT}(\varphi_z \GG) \|_{L^2(\oz)} \leq \Cs
    \sup\limits_{
    \substack{v \in H^{1}_*(\omega(z))\\\|\nabla v\|_{L^2(\oz)}=1}}
    \mathrm{Res}_z(v)
\end{align}
for the local contributions $\QQ_{z,h}$ of $\QQ_p = \sum_{z \in \mathcal{V}} \QQ_{z,h}$ from \eqref{eq:discrete_minimisation}.
The identity $\IRT(\varphi_z \GG) = \varphi_z \GG$ for $p\geq 1$ from $\varphi_z\GG\in P_{k+1}(\T(z))\subset RT_{k+p}^{{\rm pw}}(\T(z))$ allows for a $k$- and $p$-robust efficiency control of the flux reconstruction error
\begin{align}\label{eqn:alt_etaeq_eff}
    \cnst{cnst:efficiency-flux-reconstruction}^{-1} \|\QQ_p - \GG\|^2 \leq \|\nabla u-\GG\|^2 + \sum_{z \in \mathcal{V}} \osc{r}^2(\varphi_z f,\part(z))
\end{align}
with a constant $\newcnst\label{cnst:efficiency-flux-reconstruction}$ that solely depends on the shape regularity of $\mathcal{T}$.
This is deemed noteworthy and motivates two different approaches for the bound of the residual $\mathrm{Res}_z(v)$ on the right-hand side of \eqref{eq:efficiencyHdiv} for $p=0$ and $p \geq 1$ below.\medskip

\emph{Step 3.1} provides \eqref{ineq:proof-GUB-efficiency-flux} for $p=0$.
Given any normalized $v \in H^1_*(\omega(z))$ with $\|\nabla v\|_{L^2(\omega(z))}=1$,
the product $\mathcal{I}_{RT}(\varphi_z \GG) \cdot \nu_F\,v$ vanishes on any boundary facet $F \in \mathcal{F}(\partial \omega(z))$ of the patch $\omega(z)$.
This, the commuting diagram property \eqref{eq:commuting-diagram-I-RT}, and \eqref{ineq:proof-GUB-efficiency-p=0-divergence} with $\varphi_z\divpw \GG \in P_k(\mathcal{T}(z))$ in the definition of the residual $\mathrm{Res}_z(v)$ from \eqref{eq:compatibility} verify
\begin{align}\label{ineq:proof-efficiency-flux-p=0-Res}
    \mathrm{Res}_z(v) &= -(\varphi_z f + \varphi_z\divpw \GG,\Pi_k v)_{L^2(\omega(z))}\nonumber\\
    &\qquad\qquad+ \sum_{F\in\F(z) \cap \F(\Omega)} \langle  \varphi_z [\GG]_F\cdot\nu_F, \Pi_{F,k} v\rangle_{L^2(F)}.
\end{align}
The shape regularity of $\mathcal{T}$, a Poincar\'e inequality for interior vertices $z \in \mathcal{V}(\Omega)$, and a Friechrichs inequality for boundary vertices $z \in \mathcal{V}(\partial \Omega)$ provide
\begin{align}\label{ineq:proof-efficiency-GUB-p=0-scaling}
    \|h_\mathcal{T}^{-1} v\|_{L^2(\omega(z))} \approx \mathrm{diam}(\omega(z))^{-1} \|v\|_{L^2(\omega(z))} \lesssim \|\nabla v\|_{L^2(\omega(z))} = 1.
\end{align}
Given any facet $F \in \mathcal{F}(z) \cap \mathcal{F}(\Omega)$ in the facet spider $F\in\F(z)$ with facet patch $\omega(F)$,
a trace inequality thus shows \begin{align}
    \label{eqn:v_F_bound1}h_F^{-1/2}\|v\|_{L^2(F)} \lesssim \|h_\mathcal{T}^{-1} v\|_{L^2(\omega(F))}+\|\nabla v\|_{L^2(\omega(F))}\lesssim 1.
\end{align}
Cauchy inequalities on the right-hand side of \eqref{ineq:proof-efficiency-flux-p=0-Res}, the stability of the $L^2$ projection,
$\|\varphi_z\|_{L^\infty(\omega_z)} = 1$, and \eqref{ineq:proof-efficiency-GUB-p=0-scaling}--\eqref{eqn:v_F_bound1} prove that $\mathrm{Res}_z(v)$ is controlled by
\begin{align*}
    \|h_\mathcal{T}(f + \divpw \GG)\|_{L^2(\omega(z))}+
    \Big(\sum_{F \in \mathcal{F}(z) \cap \mathcal{F}(\Omega)} h_F\|[\GG]_F \cdot \nu_F\|_{L^2(F)}^2\Big)^{1/2}.
\end{align*}
The efficiency of those residual terms follows from the proof of Theorem \ref{thm:HHO_equivalence} in Section \ref{sec:residual}. In combination with \eqref{eqn:Q_q_small} and \eqref{eq:efficiencyHdiv}, this results in the global efficiency \eqref{ineq:proof-GUB-efficiency-flux} for $p = 0$ and concludes the proof of Step 3.1.\medskip

\emph{Step 3.2} provides \eqref{ineq:proof-GUB-efficiency-flux} for $p\geq 1$.
Given any normalized $v \in H^1_*(\omega(z))$ with $\|\nabla v\|_{L^2(\omega(z))}=1$, \eqref{ineq:proof-GUB-efficiency-p=0-divergence} and the identity $\IRT(\varphi_z \GG) = \varphi_z \GG$ in the definition of the residual $\mathrm{Res}_z(v)$ from \eqref{eq:compatibility} verify that $\mathrm{Res}_z(v)$ is equal to
\begin{align*}
    - (\varphi_z \mathrm{div}_\pw \GG + \Pi_{k+p}(\varphi_z f) , v)_{L^2(\omega(z))} + \sum_{T \in \mathcal{T}(z)} \langle \GG \cdot \nu_T, \varphi_z v \rangle_{L^2(\partial T)}.
\end{align*}
This, a piecewise integration by parts, and the product rule  $\nabla(\varphi_z v)=\varphi_z\nabla v+ v\nabla\varphi_z$ reveal
\begin{align}
    \mathrm{Res}_z(v) %
    &= - (\Pi_{k+p}(\varphi_z f) , v)\niLz
    + ( \GG, \nabla (\varphi_z v))\niLz.
    \label{eq:Res_z(v)}
\end{align}
The weak formulation \eqref{eq:poisson} with the test function $\varphi_z v \in H^1_0(\omega(z))\subset V$ shows $(\nabla u, \nabla(\varphi_z v))_{L^2(\omega(z))} = (f, \varphi_z v)_{L^2(\omega(z))}$. Consequently, \eqref{eq:Res_z(v)} implies
\begin{align*}
    \mathrm{Res}_z(v) =&\,
    ((1 - \Pi_{k+p})(\varphi_z f), v)\niLz
    - (\nabla u - \GG, \nabla (\varphi_z v))\niLz.
\end{align*}
This, a Cauchy-Schwarz, and a piecewise Poincar\'e inequality with the normalization $\|\nabla v\|_{L^2(\omega(z))} = 1$ provide
\begin{align}
    \label{ineq:upper-bound-Res_z(v)}
    \mathrm{Res}_z(v)  
    & \leq 
    \|\nabla u - \GG\|\niLz
    \|\nabla(\varphi_zv)\|_{L^2(\oz)}
    + C_P\osc{k+p}(\varphi_z f,\T(z)).
\end{align}
Since $\|\nabla \varphi_z\|_{L^\infty(\omega(z))}\approx h^{-1}$, the Leibniz rule, a triangle inequality, and \eqref{ineq:proof-efficiency-GUB-p=0-scaling} show that
\begin{align*}
    \|\nabla(\varphi_zv)\|_{L^2(\oz)}&\leq \|\varphi_z\|_{L^\infty(\omega(z))}\|\nabla v\|_{L^2(\omega(z))}+\|\nabla \varphi_z\|_{L^\infty(\omega(z))}\|v\|_{L^2(\omega(z))}\\ &\lesssim \|\nabla v\|_{L^2(\oz)} = 1
\end{align*}
can be bounded by a constant independent of $h$.
Thus, the combination of \eqref{eqn:Q_q_small} with \eqref{eq:efficiencyHdiv} and \eqref{ineq:upper-bound-Res_z(v)} results in the efficiency of $\|\QQ_p - \GG\|$ in \eqref{eqn:alt_etaeq_eff}
with a $k$- and $p$-robust constant \cnst{cnst:efficiency-flux-reconstruction}.
Since $\varphi_z \Pi_{k-1} f \in P_k(T) \subset P_{k+p}(T)$, the Pythagoras theorem and $\|\varphi_z\|_{L^\infty(\omega(z))} = 1$ show
$\|(1 - \Pi_{k+p})(\varphi_z f)\|_{L^2(T)} \leq \|\varphi_z(1 - \Pi_{k-1}) f\|_{L^2(T)} \leq \|(1-\Pi_{k-1})f\|_{L^2(T)}$ for all $T \in \mathcal{T}(z)$,
whence 
$$\osc{k+p}(\varphi_z f,\mathcal{T}(z)) \leq \osc{k-1}(f,\mathcal{T}(z)).$$
This and Lemma \ref{lem:local-data-oscillation} implies the efficiency 
$$\osc{k+p}(\varphi_z f, \part(z))
\lesssim \|\nablapw(u - R u_h)\|_{L^2(\omega(z))} + \osc{q}(f,\part(z))$$
of the local oscillations from \eqref{ineq:upper-bound-Res_z(v)} for any $q\in\mathbb N_0$ as in Remark \ref{rem:local_osc_eff}.
This and \eqref{eqn:alt_etaeq_eff} prove \eqref{ineq:proof-GUB-efficiency-flux} for $p \geq 1$ and conclude the proof of Step 3.2.\medskip

\emph{Step 4} affirms the efficiency of the flux reconstruction
\begin{align}\label{ineq:proof-efficiency-GUB-flux-k=0}
    \cnst{cnst:efficiency-flux-k=0}^{-1}\|\QQ_p - \GG\| \leq \trb{u - R u_h}_\pw + \osc{q}(f,\mathcal{T})
\end{align}
for the polynomial degree $k = 0$ and any $q\in\mathbb N_0$.
Let $\widetilde{u} \in V$ solve the Poisson model problem $-\Delta \widetilde{u} = \Pi_0 f$ with piecewise constant right-hand side $\Pi_0 f \in P_0(\mathcal{T})$.
A careful inspection reveals that all arguments from Step 3 apply to the case $k = 0$ for $f$ replaced by $\Pi_0 f$.
This leads to $\|\QQ_p - \GG\| \leq \cnst{cnst:efficiency-flux-k=0}\trb{\widetilde{u} - R u_h}_\pw$
with a constant $\newcnst\label{cnst:efficiency-flux-k=0}$ that solely depends on the shape of $\mathcal{T}$ (because the data oscillations on the right-hand side of \eqref{eqn:alt_etaeq_eff} vanish).
Therefore, the triangle inequality $\trb{\widetilde{u} - R u_h}_\pw \leq \trb{u - \widetilde{u}} + \trb{u - R u_h}_\pw$, the standard bound $\trb{u - \widetilde{u}} \leq C_P\osc{0}(f,\mathcal{T})$, and $C_P = 1/\pi < 1$ on simplicial domains lead to \eqref{ineq:proof-efficiency-GUB-flux-k=0}.
This and the efficiency of the data oscillations from Lemma \ref{lem:local-data-oscillation} conclude the proof of Step 4.
Notice that the constant $\cnst{cnst:efficiency-flux-k=0}$ is independent of the parameter $p$.
\medskip

\emph{Step 5} finishes the proof.
On the one hand, 
the reliability $\trb{u - R u_h}_\pw + \osc{k-1}(f,\mathcal{T}) \lesssim \eta_{\mathrm{eq},p}$ is established in Step 1--2.
On the other hand, the efficiency of the averaging operator from Theorem \ref{thm:averaging-efficiency}, Lemma \ref{lem:local-data-oscillation}, and the efficiency of the flux reconstruction in \eqref{ineq:proof-GUB-efficiency-flux} for $k \geq 1$ and in \eqref{ineq:proof-efficiency-GUB-flux-k=0} for $k = 0$ imply the efficiency
$\eta_{\mathrm{eq,p}} \lesssim \trb{u - R u_h}_\pw + \osc{q}(f,\mathcal{T})$ for any $q \in \mathbb{N}_0$.
This concludes the proof of Theorem \ref{thm:GUB-HHO}.\qed

\section{Numerical experiments}%
\label{sec:Numerical results}
This section provides numerical evidence for optimal convergence and
a comparison of the stabilization-free GUBs $\etares$ and $\eta_{\rm eq,p}$ from the Sections \ref{sec:residual} and \ref{sec:eq} with the original error estimator $\etaHHO$ from \cite[Theorem 4.3]{DiPietroDroniou2020} for the HHO method
in three 2D benchmarks. %

\subsection{A posteriori error estimation with explicit constants}
All triangulations in this section consist of right-isosceles triangles with the Poincar\'e constant $C_P = (\sqrt{2}\pi)^{-1}$ \cite{kikuchi_estimation_2007}. 
With the estimates $\Ca\leq C_\T$, $\CH\leq 1$, and $\Cb\leq C_\F$ %
from Example \ref{ex:constants}, the residual-based error estimator from \eqref{eq:residualdef} reads
\begin{align*}
    \eta^2_{\text{res}}(\T)\coloneqq& \left(C_\T\eta_{\text{res}, 1}(\T) + C_P\eta_{\text{res}, 2}(\T)
+C_\F\eta_{\text{res}, 3}(\T)\right)^2+ C_\F^2\eta^2_{\text{res}, 4}(\T).
\end{align*}
The equilibrated GUB from Theorem \ref{thm:GUB-HHO},
\begin{align*}
	\eta_{\mathrm{eq},p}^2(\T)\coloneqq \left(C_P\osc{k+p}(f, \T) + \| Q_p^\Delta \|_{L^2(\Omega)}\right)^2 + \trb{(1-\mathcal{A}) Ru_h}_\pw^2,
\end{align*}
depends on the post-processed quantity $Q_p^\Delta\coloneqq Q_p-\nablapw Ru_h$ with the Raviart-Thomas function $Q_p\in RT_{k+p}(\T)$ of degree $k+p$ for $p\in\mathbb N_0$ from Subsection \ref{sub:Computation_Q}.
An algorithmic description of the computation of $Q_p^\Delta$ for arbitrary $p$ follows in Appendix A.
Theorem \ref{thm:GUB-HHO} shows that $\etaeqnew$ is efficient and reliable for all $p\in\mathbb N_0$.
The residual-based error estimator $\eta_{\rm HHO}$ from \cite[Theorem 4.3]{DiPietroDroniou2020} reads
\begin{align*}
	\eta_{\rm HHO}^2(\T) \coloneqq& \sum_{T\in\T}  \bigg( C_Ph_{T}\| (I - \Pi_{T,0})(f+{\Delta}_{\rm pw} Ru_h)\|_{L^2(T)}\\
	&\quad+ \sqrt{\sum_{F\in\F(T)} C_{\partial T} h_T \|R_{T, F}^k u_h\|_{L^2(F)}^2}\bigg)^2
	 + \trb{(I-\mathcal{A}) Ru_h}_\pw^2
\end{align*}
with the operator $R_{T, F}^k$ from \cite[Eq.~(2.59)]{DiPietroDroniou2020} for the original HHO stabilization \cite[Eq.~(2.22)]{DiPietroDroniou2020} that induces the global stabilization $s_h$.
The constant %
$C_{\partial T}=12C_P(C_P+\Ctr)\leq 2.0315$ for right-isosceles triangles bounds the trace of $f-\Pi_0 f$ for $f\in H^1(T)$ and improves on the estimate $C_{\partial T}\leq (C_P(h_T|\partial T|/|T|)(1+C_P))^2=2.524$ from \cite[Subs.~4.1.1]{DiPietroDroniou2020}.
\begin{lemma}[Poincar\'e-type inequality on trace]
Given a simplex $T \subset \mathbb{R}^n$, any $f\in H^1(T)$ and $C_{\partial T}\coloneqq (n+1)h_T^2|T|^{-1} C_P(C_P + 2\Ctr/n)$ satisfy
\begin{align}\label{eqn:C_pT}
    \|f-\Pi_{0} f\|_{L^2(\partial T)}^2 \leq C_{\partial T}h_T\|\nabla f\|_{L^2(T)}^2.
\end{align}
\end{lemma}
\begin{proof}
Abbreviate $\tilde \ell(F)\coloneqq (n+1)h_Fh_T^2|T|^{-1}$ for the facet $F\in\F(T)$ of $T$ and apply Lemma \ref{lem:Skeletal trace inequality} to the singleton triangulation $\{T\}$ and $f-\Pi_0f\in H^1(T)$.
This and the Poincar\'e inequality $\|f-\Pi_{0, T}f\|_{L^2(T)}\leq C_P h_T\|\nabla f\|_{L^2(T)}$ show
\begin{align*}
 \sum_{F\in\F(T)} \tilde \ell(F)^{-1}\|f-\Pi_{0, T} f\|_{L^2(F)}^2
 &\leq(C_P+2\Ctr/n)C_P \|\nabla f\|_{L^2(T)}^2.
\end{align*}
The assertion \eqref{eqn:C_pT} follows from the observation that $\tilde \ell(F)$ is maximized on the facet $F\in\F(T)$ with $h_F=h_T$.\qed
\end{proof}

\subsection{Implementation and adaptive algorithm} %
\label{sub:Implementation}

Our implementation of the HHO method in MATLAB  uses nodal bases for the spaces $P_k(\F), P_k(\T)$, and $P_{k+1}(\T)$ and the direct solver \textit{mldivide} (behind the \textbackslash-operator) for the discrete system of equations representing \eqref{eqn:HHO}.
For implementation details on the HHO method itself we refer to  \cite[Appendix B]{DiPietroDroniou2020}.
The integration of polynomial expressions is carried out exactly.
%
The errors in approximating non-polynomial expressions, such as exact solutions $u$ and source terms $f$, by polynomials of sufficiently high degree are expected to be very small and are neglected for simplicity.

Algorithm \ref{alg:afem} displays the standard adaptive algorithm (AFEM)  \cite{carstensen_axioms_2014,carstensen_axioms_2017} %
driven by the refinement indicators, for any triangle $T\in \T$,
\begin{align}\label{eqn:etares_T}
    \eta^2_{\rm res}(T)&\coloneqq |T|\|f+{\Delta}_{\rm pw} Ru_h\|^2_{L^2(T)}
+|T|^{1/2}\sum_{F\in \cF(T)}\|[ {\nabla}_{\rm pw} Ru_h ]_F \|_{L^2(F)}^2 %
\end{align}
with the modified jump %
$[\bullet]_F\coloneqq \bullet\times n_F$ along a boundary side $F\in\F(\partial \Omega)$ and the newest-vertex-bisection (NVB).
The sum of this over all triangles is, up to some multiplicative constants, equivalent to $\eta_{\rm
res}^2(\mathcal{T})$.
\newcommand{\MM}{\mathcal M}
\begin{algorithm}
\caption{AFEM algorithm}
\label{alg:afem}
\textbf{Input:} Initial regular triangulation $\T_0$ and polynomial degree $k\in\mathbb N_0$ of the HHO method
\begin{algorithmic}
    \For{levels $\ell\coloneqq 0, 1, 2...$}
    \State \textbf{Solve} \eqref{eqn:HHO} for discrete solution $u_\ell\in V_\ell$ exactly on $\T_\ell$ and compute $Ru_\ell$
    \State \textbf{Compute} (refinement indicators) $\etares^2(T)$ for all $T\in\T_\ell$ %
    \State \textbf{Mark} minimal subset $\MM_\ell\subset\T_\ell$ with $\frac12\sum_{T\in\T_\ell}\etares^2(T) \leq \sum_{T\in\MM_\ell}\etares^2(T)$
    \State \textbf{Refine} $\T_\ell$ to smallest NVB refinement $\T_{\ell+1}$ with $\MM_\ell\subseteq\T_\ell\setminus\T_{\ell+1}$
    \EndFor
\end{algorithmic}
 \textbf{Output:} sequences of triangulations $\T_\ell$ and $Ru_\ell$
\end{algorithm}

\subsection{High oscillations on the unit square}%
\begin{figure}
		\centering
		\includegraphics{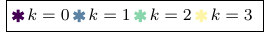}\\
		\hspace*{-2em}\hbox{\includegraphics{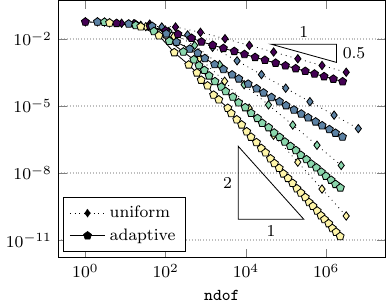}
		\includegraphics{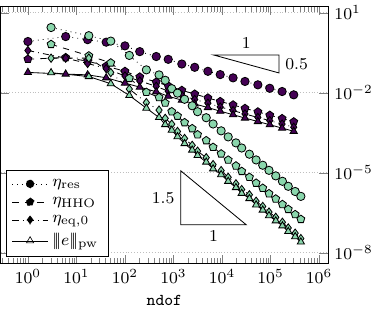}}
		\caption{Convergence history of the energy error $\trb{e}_{\rm pw}$ (left) and the GUBs $\etares, \etaHHO, \etaeqnewa{0}$ (right) on the square domain}
		\label{fig:SquareConvergence}
	\end{figure}

\label{sub:Smooth solution}
This benchmark on the unit square $\Omega\coloneqq(0,1)^2$ considers the Laplace equation $-\Delta u=f$ with source term $f$ matching the smooth exact solution
$$u(x, y) = x(x-1)y(y-1) e^{-100((x-1/2)^2 + (y-117/1000)^2)}\in C^\infty(\Omega).$$
Figure \ref{fig:SquareConvergence} displays the energy norm $\trb{e}_\pw$ of the error $e\coloneqq u-Ru_h$ for uniform and adaptive refinement by Algorithm \ref{alg:afem} on the left.
The smooth solution allows for optimal convergence rates $(k+1)/2$ in the number \ndof\ of degrees of freedom, while the adaptive mesh sequence leads to a lower energy error with respect to \ndof.
The GUBs $\etares, \etaeqnewa p$, and $\etaHHO$ are efficient and therefore equivalent to the energy error $\trb{e}_\pw$, see Figure \ref{fig:SquareConvergence} on the right for $k=0, 2$.

\begin{figure}
		\centering
		\includegraphics{./Figures/Legend_k_short_out.pdf}\\
		\includegraphics{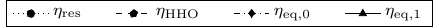}\\
		\hbox{
        \includegraphics{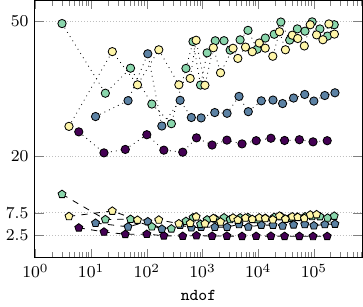}
        \includegraphics{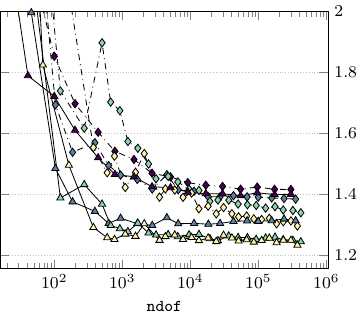}}
		\caption{History of the overestimation factor
		$EF(\eta)=\eta/\trb{e}_\pw$ for the residual-based error estimators $\etares, \etaHHO$ (left) and $\etaeqnewa0, \etaeqnewa1$ (right) on the unit square}
		\label{fig:SquareGUB}
	\end{figure}

Figure \ref{fig:SquareGUB} shows the efficiency indices $EF(\eta)\coloneqq \eta/\trb{e}_\pw$ for the residual-based GUBs $\eta=\etares, \etaHHO$ and the equilibration-based GUB $\eta=\etaeqnewa p$ for $p=0,1$. Higher values of $p$ for a more expensive postprocessing in $\etaeqnew$ do not significantly improve on $\etaeqnewa1$.

\subsection{Analytical solution for the slit domain}%
\begin{figure}
		\centering
		\includegraphics{./Figures/Legend_k_short_out.pdf}\\
		\hspace*{-2em}\hbox{\includegraphics{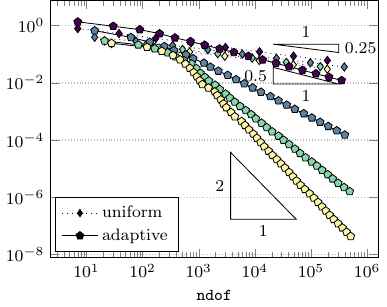}
		\includegraphics{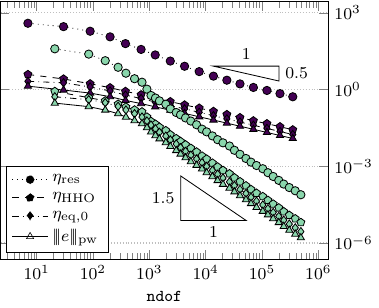}}
		\caption{Convergence history of the energy error $\trb{e}_{\rm pw}$ (left) and the GUBs $\etares, \etaHHO, \etaeqnewa{0}$ (right) on the slit domain}
		\label{fig:SlitConvergence}
	\end{figure}

\label{sub:Slit singularity on the slit domain}
The source term $f\in L^2(\Omega)$ in the second benchmark on the slit domain $\Omega\coloneqq(0,1)^2\setminus([0,1)\times\{0\})$
matches the singular solution (in polar coordinates)
$$u(r, \varphi) = r^{1/2}(r^2\sin(\varphi)^2 - 1)(r^2\cos(\varphi)^2 - 1)\sin(\varphi/2).$$
The singularity of $u$ at the origin $(0,0)$ leads to reduced convergence rates $1/4$ under uniform refinement, regardless of the polynomial degree $k$.
Figure \ref{fig:SlitConvergence} shows that the adaptive algorithm recovers optimal rates and verifies the equivalence of the GUBs $\etares, \etaeqnewa p$, and $\etaHHO$ to the energy error $\trb{e}_\pw$.

\begin{figure}
		\centering
		\includegraphics{./Figures/Legend_k_short_out.pdf}\\
		\includegraphics{./Figures/Legend_ErrorV2_out.pdf}\\
		\hbox{
        \includegraphics{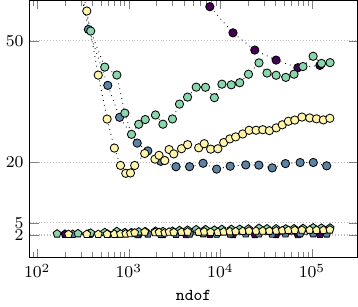}
        \includegraphics{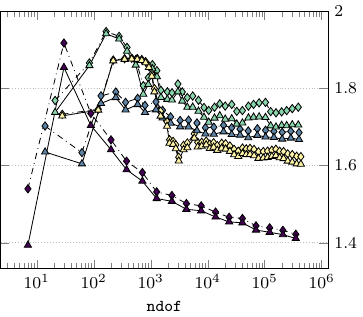}}
		\caption{History of the overestimation factor
		$EF(\eta)=\eta/\trb{e}_\pw$ for the residual-based error estimators $\etares, \etaHHO$ (left) and $\etaeqnewa0, \etaeqnewa1$ (right) on the slit domain}
		\label{fig:SlitGUB}
	\end{figure}

The efficiency indices in Figure \ref{fig:SlitGUB} show a strong overestimaton by $\etares$ in the preasymptotic regime (undisplayed) with values $EF(\etares)>60$.
However, asymptotically the quotients $EF(\eta)=\eta/\trb{e}_\pw$ for the two residual-based GUBs $\eta_{\rm res}, \eta_{\rm HHO}$ differ only by a factor $10$, while the equilibrated GUBs $\etaeqnew$ provide the closest values to $1$.
\subsection{Corner singularity in the L-shaped domain}%
\begin{figure}
		\centering
		\includegraphics{./Figures/Legend_k_short_out.pdf}\\
		\hspace*{-2em}\hbox{\includegraphics{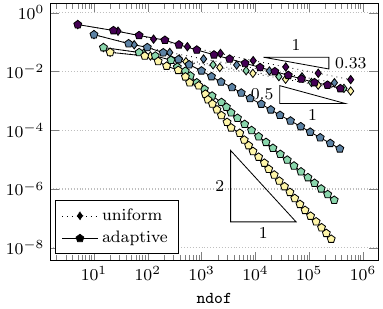}
		\includegraphics{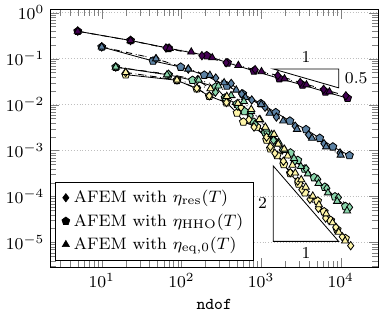}}
		\caption{Convergence history plot of the energy error $\trb{e}_\pw$ on the L-shaped domain with uniform and adaptive refinement with AFEM, driven by $\etares(T)$, (left) and for AFEM, driven by $\etares(T), \etaHHO(T)$, and $\etaeqnewa0(T)$, (right).}
		\label{fig:LshapeConvergence}
	\end{figure}

\begin{figure}
		\centering
		\hbox{%
		\includegraphics{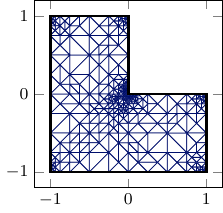}
		\includegraphics{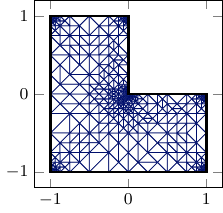}
		\includegraphics{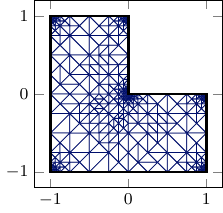}}

		\caption{Adaptive triangulations on the L-shaped domain for $k=3$ from AFEM, driven by $\etares(T)$ (left, $|\T|=882$), driven by $\etaHHO(T)$ (middle, $|\T|=907$), and driven by $\etaeqnewa0(T)$ (right, $|\T|=919$).}
		\label{fig:LshapeMesh}
	\end{figure}

\label{sub:Singular solution}
The third benchmark problem is set in the 
L-shaped domain $\Omega = (-1, 1)^2\setminus [0, 1)^2$ with constant right-hand side $f\equiv1$
with an exact solution $u \in H^{1+s}(\Omega)$ for all $0 \leq s < 2/3$ \cite[Theorem 14.6]{Dauge1988}.
Figure \ref{fig:LshapeConvergence} 
displays the convergence history of the error $e\coloneqq u-Ru_h$ and compares the adaptive scheme, Algorithm \ref{alg:afem}, driven by the refinement indicators $\eta_{\rm res}(T)$ from \eqref{eqn:etares_T} and
\begin{align*}
    \etaHHO^2(T)&\coloneqq |T|\| (I - \Pi_0)(f+{\Delta}_{\rm pw} Ru_h)\|_{L^2(T)}^2 + \|\nabla(1-\mathcal{A}) Ru_h\|_{L^2(T)}^2\notag \\
    &\quad+ 
	|T|^{1/2}\sum_{F\in\F(T)} \| R_{T, F}^k u_h\|_{L^2(F)}^2,\\
	\etaeqnewa0^2(T)&\coloneqq\osc{k}^2(f, T) + \| Q_0^\Delta \|_{L^2(T)}^2 + \|\nabla(1-\mathcal{A}) Ru_h\|_{L^2(T)}^2
\end{align*}
for $T\in\T$ that are induced from the GUB $\etares, \etaHHO$, and $\etaeqnewa0$.
Here, the norm $\trb{e}_\pw$ of the distance $e$ from the discrete solution $Ru_h\in P_{k+1}(\T)$ over $\T$ to the unknown solution $u\in H^1(\Omega)$ is approximated by $\trb{\widehat u-Ru_h}_\pw$, where $\widehat u\in P_{k+1}(\widehat \T)$ is the HHO approximation of $u$ on an adaptive refinement $\widehat \T$ of $\T$ with at least $2|\T|\leq|\widehat\T|$ elements.
The three adaptive schemes (Algorithm \ref{alg:afem}, driven by $\etares(T), \etaHHO(T)$, or $\etaeqnewa0(T)$) recover optimal rates of convergence and lead to similar local refinement of the adaptive mesh sequences as in Figure~\ref{fig:LshapeMesh}.

\subsection{Conclusion}
The adaptive mesh-refining algorithm recovers optimal convergence rates in all three benchmarks.
This holds for
AFEM driven by any of the three refinement indicators derived from the GUB $\etares, \etaHHO$, and $\etaeqnew$.
The generated mesh sequences from the adaptive schemes, driven by the different estimators, display a very similar concentration of the local mesh-refinement as in Figure~\ref{fig:LshapeMesh}.
All three benchmarks verify that the considered error estimators are GUB with reliability constant $1$, while the post-processing in the equilibrated GUB $\etaeqnew$ produces minimal overestimation with significant additional computational costs.
%

%
%
%
%
%
%
%
%
%
%
%
%
%
%
%

%
%
%
%
%
%
%

%
%
%
%
%
%
%
%
%
%
%
%
%
%
%
%
%
%
%
%
%
%
%
%
%
%
%
%
%
%
%
%
%
%
%
%
%
%
%
%
%
%
%
%
%
%

%
%
%
%
%
%
%
%
%
%
%
%
%
%
%
%
%
%
%

\bibliographystyle{spmpsci}      
\bibliography{Bibliography}   


\newpage
\appendix %
\renewcommand{\thesection}{Appendix \Alph{section}:}
\section{Equilibration algorithm for higher order}
\renewcommand{\thesection}{\Alph{section}}
\label{app:equilibration}
\FloatBarrier
\newcommand{\inda}{a}
\newcommand{\kp}{{q}}
\newcommand{\tszD}{\tilde\sigma_z^\Delta}
\newcommand{\szD}{\sigma_z^\Delta}
\newcommand{\sz}{\sigma_z}
\newcommand{\RTkz}{RT^{0}_k}
\newcommand{\RTkpwz}{RT_{\kp}^{{\rm pw}, 0}}
\newcommand{\RTpwk}{RT_{\kp}^{{\rm pw}}}
\renewcommand{\PiRT}{\Pi_{RT}}
\newcommand{\dofBasis}{\Lambda_T}

\newcommand{\dx}{\mathrm{dx}}
\newcommand{\ds}{\mathrm{ds}}
\newcommand{\Llm}{\lambda_{T, \alpha}^{\ell, m}}
\newcommand{\Ldivlm}{\lambda_{T, \ddiv}^{\ell, m}}
\newcommand{\LFj}{\lambda_{T, E}^j}
\newcommand{\Lrestr}{\lambda_T^{r}}
\newcommand{\tBRT}{\widetilde{\dofBasis}}
\newcommand{\BRT}{{\mathcal{B}}_{RT}}
\newcommand{\BRTT}{{\mathcal{B}}_{RT,T}}

\newcommand{\cFj}{c_{T, E}^j}
\newcommand{\cdivlm}{c_{T, \ddiv}^{\ell, m}}
\newcommand{\cdivlma}{c_{T_\inda, \ddiv}^{\ell, m}}
\newcommand{\crestr}{c_{T}^r}
\newcommand{\vfFj}{\varphi_{T, E}^j}
\newcommand{\vfdivlm}{\varphi_{T, \ddiv}^{\ell, m}}
\newcommand{\vfrestr}{\varphi_{T}^r}
\newcommand{\Sz}{\mathcal S(z)}
\newcommand{\Vz}{V(z)}
\newcommand{\szd}{\sigma_z^\Delta}

The post-processed quantity $\QQ_p\in RT_{k+p}(\T)$ from Subsection \ref{sub:Computation_Q} enters the equilibrated error estimator $\etaeqnew(\T)$ in Theorem \ref{thm:GUB-HHO} and could be computed by a minimization problem on the vertex patches.
The solution property \eqref{eqn:HHO_solution_property} gives rise to the two cases $r = 0$ if $k = 0$ and $r = k+p$ if $k \geq 1$ for the polynomial degree $r$ in the equilibrium $\ddiv \QQ_p + \Pi_{r} f = 0$ in $\Omega$ from \eqref{eqn:div_Q_f}.
This appendix follows \cite{verfurth_note_2009,cai_robust_2012,bertrand_weakly_2019,braess_finite_2007} to compute the quantity of interest
$Q_p - \nabla_\pw R u_h$ directly in an efficient two-step procedure in 2D.
Throughout this appendix, fix $k,p\in \mathbb N_0$ and abbreviate $\kp\coloneqq k+p$ and $\GG \coloneqq \nabla_\pw R u_h\in P_{k}(\T;\R^2)$.
Let the data $f\in L^2(\Omega)$ be given and assume, for the sake of brevity, that $f \in P_0(\mathcal{T})$ if $k = 0$.

\subsection{Overview}
\begin{figure}
        \centering
        \scalebox{.6}{\includegraphics{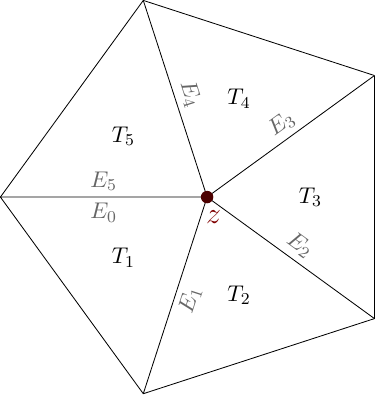}}
        \scalebox{.6}{\includegraphics{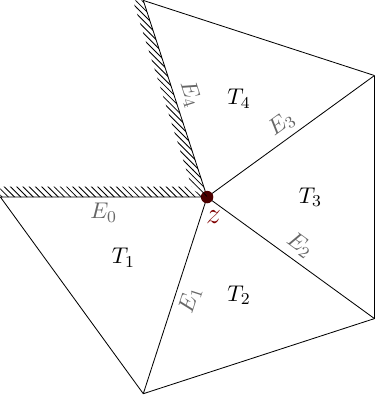}}
        \caption{Triangulation $\T(z)$ and enumeration of the edges $\F(z)$ of the vertex-patch $\oz$ for an internal vertex $z\in\mathcal V(\Omega)$ with $N=|\T(z)|=5$ (left) or boundary vertex $z\in\mathcal V(\partial \Omega)$ with $N=|\T(z)|=4$ (right).}
        \label{fig:NodePatch}
\end{figure}

Recall the definition \eqref{eq:discrete_minimisation} of the summand $\Rzh$ in $Q_p \coloneqq\sum_{z\in \mathcal V}\Rzh$ from Subsection \ref{sub:Computation_Q}
with the piecewise Raviart-Thomas interpolation $\IRT:H^1(\T;\R^2)\to \RTpwk(\T)$ \cite[Section III.3.1]{boffi_mixed_2013}.
The focus is on one vertex $z\in\mathcal V$ with vertex-patch $\omega(z)$ and its triangulation $\T(z)=\{T\in\T\ |\ z\in T\}$.
Consider set of edges $\mathcal{F}$ and the facet-spider $\F(z) =\{E\in\F\ |\ z\in E\}$ as in Figure \ref{fig:NodePatch}. The nodal basis function $\varphi_z\in S^1(\T(z))$ gives rise to the discrete spaces
\begin{align*}
    &\RTkpwz(\mathcal{T}(z))\coloneqq\{\sz\in RT_{\kp}^{\text{pw}}(\T(z))\ :\ \sigma_{z}\cdot \nu_E=0 \text{ for all }E\in\F\setminus\F(z)\},\\
    &\Sz\coloneqq\left\{\sz \in \RTkpwz(z)\ :\ \begin{array}{llr}
    \ddiv\sz &=-\Pi_{T,{\kp}}(\varphi_z(f+\ddiv \GG))&\text{ for all } T\in\T(z)\\
    {[}\sz\cdot \nu_E{]_E} &= -\Pi_{E, {\kp}}(\varphi_z[\GG \cdot \nu_E]_E)&\text{ for all } E\in\F(\omega(z))
    \end{array}\right\}.
\end{align*}

\begin{proposition}[alternative minimization]\label{prop:RzhD_def}
It holds
\begin{align}\label{eqn:RzhD_def}
   \RzhD\coloneqq\Rzh - \IRT(\varphi_z \GG) = \argmin\limits_{%
\sz \in\Sz
     }
 \| \sz\|\niLz.
\end{align}
\end{proposition}
\begin{proof}
Recall $f_z = \Pi_{\tilde p}(\varphi_z f - \GG \cdot \nabla \varphi_z)$ from \eqref{def:f-z}. (Notice that this formula coincide with the definition \eqref{def:f-z} for $k = 0$ because $f \in P_0(\mathcal{T})$.) Given any $\sigma_z \in \mathcal{S}(z)$, the commuting diagram property $\divpw\circ \IRT = \Pi_{\kp}\circ\divpw$ \cite[Proposition 2.5.2]{boffi_mixed_2013} shows
\begin{align}\label{eq:divergence-S(z)-IRT}
    \mathrm{div}_\pw (\sigma_z + \mathcal{I}_{RT} (\varphi_z \GG)) = \mathrm{div}_\pw \sigma_z + \Pi_{\tilde p}\,\mathrm{div}_\pw (\varphi_z \GG) = - f_z.
\end{align}
By design of the interpolation $\mathcal{I}_{RT}$, $(\mathcal{I}_{RT} (\varphi_z \GG)_{|T}\cdot\nu_E)_{|E} = \Pi_{E,\kp}(\varphi_z \GG)_{|T}\cdot\nu_E$ holds and so
$[\IRT (\varphi_z \GG)\cdot \nu_E]_E = \Pi_{E, {\kp}}(\varphi_z[\GG \cdot \nu_E]_E)$ follows for any $E \in \mathcal{F}(T)$ and $T\in\T$.
Therefore, the jump $[\sigma_z + \IRT(\varphi_z \GG)]_E\cdot\nu_E \equiv 0$ vanishes on $E \in \mathcal{F}(\omega(z))$, whence $\sigma_z + \IRT(\varphi_z \GG) \in RT_{\tilde{p}}(\mathcal{T}(z))$.
Since $RT_{\tilde{p}}^{\pw,0}(\mathcal{T}(z)) \cap H(\mathrm{div},\omega(z)) = RT_{\tilde{p}}^0(\mathcal{T}(z))$, this and \eqref{eq:divergence-S(z)-IRT} imply $\sigma_z + \IRT(\varphi_z \GG) \in \mathcal{Q}_h(z)$ for any $\sigma_z \in \mathcal{S}(z)$.
In particular, $\Sz + \IRT(\varphi_z \GG) \subseteq \mathcal{Q}_h(z)$.
On the other hand, similar arguments verify the reverse inclusion $\mathcal{Q}_h(z) \subseteq \Sz + \IRT(\varphi_z \GG)$. %
The substitution $\Sz + \IRT(\varphi_z \GG) = \mathcal{Q}_h(z)$ in \eqref{eq:discrete_minimisation} concludes the proof.\qed
\end{proof}
This establishes that the norm $\|Q^\Delta_p\|$ of $Q^\Delta_p\coloneqq\sum_{z\in \mathcal V}\RzhD = Q_p - \GG$ contributes to the equilibrated error estimator and the remaining parts of this
 appendix compute the minimizer $\RzhD$ of \eqref{eqn:RzhD_def} in a two-step procedure.
First, Algorithm \ref{alg:eq} generates the coefficients of a particular solution $\tszD\in \Sz$ in terms of the finite element basis $\BRT$ of $\RTpwk(\mathcal{T}(z))$ from Subsection \ref{sub:basisRT}.
The second step computes the correction \begin{align}\label{eqn:linear_min}
    \szD\coloneqq\RzhD - \tszD= \argmin_{\sz\in \Vz}\|\sz +\tszD \|_{L^2(\omega(z))}
\end{align}
in terms of the low-dimensional unconstrained
minimization problem 
over the linear space $\Vz\coloneqq\Sz-\tszD=\RTpwk(\mathcal{T}(z))\cap H(\ddiv=0, \omega(z))$ characterized in Lemma \ref{lem:appendixB}.
Because \cite[Lemma 3.1]{cai_robust_2012} is wrong (take, e.g.,  $\tau_K=\mathrm{curl}\,b_K\ne 0$ for the element bubble function $b_K$ in their notation to see that uniqueness for general polynomial degrees $\kp$ cannot hold) and \cite{cai_robust_2012} omits algorithmic details, this appendix focuses on the explicit characterization of the degrees of freedom for the minimization problem \eqref{eqn:linear_min} over $\Vz$.
\subsection{Degrees of freedom for $\RTpwk(T)$}
\label{sub:basisRT}
This subsection introduces a basis for the Raviart-Thomas finite element on $T\in\T(z)$ and starts with the definition of some linear functionals on $H(\ddiv, \T)$. For any $\sigma\in H(\ddiv, \T)$, set
\begin{align*}
    \Llm(\sigma) \coloneqq& \int_T  \sigma\cdot e_\alpha\ x_1^\ell x_2^m \dx&& 0\leq \ell+m\leq \kp-1,\alpha=1, 2,\\
    \Ldivlm(\sigma) \coloneqq& \int_T \ddiv \sigma \ x_1^\ell x_2^m \dx&& 0\leq \ell+m\leq \kp,\\
    \LFj(\sigma) \coloneqq& \int_E \sigma_{|T} \cdot \nu_T\ s^j \ds&& 0\leq j\leq \kp, E\in\F(T).
\end{align*}
Here and throughout, $e_1=(1,0)$ and $e_2=(0,1)$ denote the canonical unit vectors in $\R^2$.
Note that the (classical) degrees of freedom for the Raviart-Thomas finite element $RT_\kp(T)$ of degree $\kp\in\mathbb N_0$ from \cite{brezzi_mixed_1991} read
$$\tBRT\coloneqq\{\LFj, \Llm\ :\ \text{for }0\leq j\leq \kp, 0\leq \ell+m\leq \kp-1, \alpha=1,2, E\in\F(T)\}.$$
This appendix requires, for the construction of $\tszD\in\Sz$,  a different set of (unisolvent) degrees of freedom $\dofBasis$ for $RT_\kp(T)$ that includes the edge and divergence moments 
\begin{align}\label{eqn:dof_basis}
    \dofBasis^0\coloneqq \{\LFj, \Ldivlm\ :\ \text{for }0\leq j\leq \kp, 1\leq \ell+m\leq \kp, E\in\F(T)\}\subseteq \dofBasis.
\end{align}
(The set $\dofBasis^0$ itself is linear independent \cite[Lemma 3.1]{verfurth_note_2009}.) 
Given any $\dofBasis$ with \eqref{eqn:dof_basis}, denote the remaining $N_\kp=\kp(\kp-1)/2$ degrees of freedom $\dofBasis\setminus\dofBasis^0$ by $\Lrestr$ for $r=1,\dots,N_\kp$.
Let $\BRTT=\{\vfFj, \vfdivlm, \vfrestr\}$ be the unique basis of $RT_\kp(T)$ dual to $\dofBasis$ with inferred indices from $\dofBasis$.
Then, the collection $\BRT\coloneqq \bigcup_{T\in\T(z)} \BRTT$ is a basis of $\RTpwk(\mathcal{T}(z))$ and any function $\sz\in \RTpwk(\mathcal{T}(z))$ has the representation 
\begin{align}\label{eqn:basis_RTkpwz}
 \sz \coloneqq \sum_{T\in\T(z)}\left(\sum_{E\in\F(T)}\sum_{j=0}^\kp \cFj\vfFj + \sum_{1\leq \ell+m\leq \kp}\cdivlm\vfdivlm + \sum_{r=1}^{N_\kp}\crestr\vfrestr\right)
 \end{align}
 with coefficients $\cFj = \LFj(\sz),  \cdivlm = \Ldivlm(\sz)$, and $\crestr = \Lrestr(\sz)$ for all $T\in\T(z), E\in\F(T)$, and $0\leq j\leq \kp, 1\leq \ell+m\leq \kp, 1\leq r\leq N_\kp$.
 By duality, the coefficients $c_{T, E}^j$ with $0\leq j\leq \kp$ uniquely determine the normal trace $(\sigma_{z|T})_{|E}\cdot \nu_E\in P_\kp(E)$ on the edge $E\in\F(T)$ of $T\in\T(z)$.
 Any set of degrees of freedom $\dofBasis$ with \eqref{eqn:dof_basis} works with the equilibration algorithm in \ref{sub:app:alg}.
\newline
\begin{example}[Construction of $\dofBasis$]
This example presents a generic procedure to obtain such a set from $\tBRT$.
The integration by parts formula shows that the lowest-order divergence moment $\lambda_{T, \div}^{0,0} = \sum_{E\in\F(T)} \lambda_{T, E}^0$  depends linearly on the lowest-order edge moments and, similarly, the sums
\begin{align}
    \label{eqn:ibp}
    \Ldivlm + \ell\lambda_{T, 1}^{\ell-1, m} + m\lambda_{T, 2}^{\ell, m-1}\in (P_{\ell + m}(\F(T)))^*
\end{align}
are functionals on $P_{\ell + m}(\F(T))$ (summands with negative indices are understood as zero).
This relation allows for the substitution of volume moments in $\tBRT$ for divergence moments $\Ldivlm$, $1\leq \ell+m\leq \kp$, and leads to $\dofBasis$ with \eqref{eqn:dof_basis}. 
The remaining degrees of freedom $\dofBasis\setminus\dofBasis^0$ are volume moments of the form $\Lrestr=\Llm$ for a fixed $\alpha \in \{1,2\}$, e.g.,
$$\dofBasis \setminus \dofBasis^0=\{\lambda_{T, 2}^{\ell, m}\ :\ 1\leq \ell\leq \kp-1, 0\leq m\leq \kp-1-\ell\}.$$
\end{example}

\subsection{Equilibration algorithm}\label{sub:app:alg}
\FloatBarrier
\begin{algorithm}
\caption{Particular solution in $\Sz$}
\label{alg:eq}
\textbf{Input:} Data $f\in L^2(\oz)$ and $\GG \in H^1(\T(z))^2$ for vertex $z\in\mathcal V$.
\begin{algorithmic}[1]
\State Initialize all coefficients $\cFj,\cdivlm, \crestr$ in \eqref{eqn:basis_RTkpwz} with zero.
    \For{$\inda\coloneqq 1:N$}
    \State $c_{T_\inda, E_{\inda-1}}^{0}\coloneqq \begin{cases}
    0%
    &\text{if } \inda=1,\\
    ([\GG]_E\cdot \nu_E, \varphi_z)_{L^2(E_{\inda-1})} - c_{T_{\inda-1}, E_{\inda-1}}^{0} &\text{else}
    \end{cases}$
    \State $c_{T_\inda, E_{\inda\phantom {-1}}}^{0}\coloneqq (f+\ddiv \GG, \varphi_z)_{L^2(T_\inda)} - c_{T_\inda, E_{\inda-1}}^{0}$
    \For{$1\leq \ell+m\leq \kp$}
        \State $c_{T_\inda, \div}^{\ell,m}\coloneqq (f+\ddiv \GG, \varphi_zx_1^\ell x_2^m)_{L^2(T_\inda)}$
    \EndFor
    \For{$1\leq j\leq \kp$}
        \State $c_{T_\inda, E_{\inda-1}}^{j}\coloneqq 0 %
        $
    \State $c_{T_\inda, E_{\inda\phantom {-1}}}^{j}\coloneqq ([\GG]_E\cdot \nu_E, \varphi_zs^j)_{L^2(E_\inda)}%
    $
    \EndFor
    \EndFor
\end{algorithmic}
\textbf{Output:} $\tszD\in\RTpwk(\mathcal{T}(z))$ defined %
by \eqref{eqn:basis_RTkpwz} with coefficients $\cFj,\cdivlm, \crestr$.
\end{algorithm}
This subsection presents the equilibration procedure, starting with Algorithm \ref{alg:eq}, that computes an admissible function $\tszD\in\Sz$ in terms of the representation \eqref{eqn:basis_RTkpwz}.
Enumerate the $N\coloneqq|\T(z)|$ triangles $T\in\T(z)$ from $1$ to $N$ as in  Figure \ref{fig:NodePatch}.
Any two neighbouring triangles $T_\inda, T_{\inda+1}$ share an edge $E_\inda\coloneqq T_\inda\cap T_{\inda+1}$ for $\inda=1,...,N-1$.
If $z\in\mathcal V(\Omega)$ is an interior vertex, $T_1$ and $T_{N}$ share an additional edge $E_0\coloneqq E_{N}\coloneqq T_1\cap T_{N}$.
For a boundary vertex $z\in\mathcal V(\partial\Omega)$, $T_1, T_{N}$ have the distinct boundary edges $E_0, E_{N}\in \F(z)\cap \F(\partial\Omega)$.
The following lemma shows correctness of Algorithm \ref{alg:eq} under the compatibility condition \eqref{eq:solution_property} and represents step one of the equilibration algorithm.
The final step is the local minimization problem in Lemma \ref{lem:appendixB} that provides $\RzhD$ from \eqref{eqn:RzhD_def}.
Both proofs are provided in \ref{sub:app:proofs}.

\begin{lemma}\label{lem:appendix}
Given $z\in\mathcal V$, let $\{\vfFj, \vfdivlm, \vfrestr\}$ be the basis of $RT_\kp(T)$ dual to $\dofBasis$ with \eqref{eqn:dof_basis} for all $T\in\T(z)$.
Suppose $f \in L^2(\oz)$ and $\GG\in P_\kp(\T(z);\mathbb{R}^2)$ satisfy %
\begin{equation}\label{eqn:crutial_identity}
(\GG, \nabla \varphi_z)_{L^2(\oz)}=(f,\varphi_z)_{L^2(\oz)}\quad\text{if }z\in\mathcal V(\Omega).
\end{equation}
 Then the output of Algorithm \ref{alg:eq} with input $f$ and $\GG$ defines a function $\tszD\in\Sz$.
 \end{lemma}
 Note that \eqref{eqn:crutial_identity} is a local version of \eqref{eq:solution_property} and therefore holds for the HHO method with the choice $\GG \coloneqq \nabla_\pw R u_h$ as proven in \eqref{eqn:HHO_solution_property}.
 This allows for the computation of $Q^\Delta_p \coloneqq\sum_{z\in \mathcal V}\RzhD = Q_p - \GG$ in terms of local and unconstrained minimization problems on the vertex-patches $\oz$.
 \begin{lemma}\label{lem:appendixB}
Given $z\in\mathcal V$, let $\{\vfFj, \vfdivlm, \vfrestr\}$ be as in Lemma \ref{lem:appendix} for all $T\in\T(z)$ and let $\tszD\in \Sz$ be arbitrary.
 Then $\Vz\coloneqq \Sz - \tszD$ is a linear vector space and consists of all functions of the form
 \begin{align} \label{eqn:sigma}
 \sum_{\inda=1}^{N}\left(d_0 (\varphi_{T_\inda, E_{\inda-1}}^0-\varphi_{T_{\inda}, E_{\inda}}^0)
 + \sum_{\ell=1}^\kp(d_{E_{\inda-1}}^\ell\varphi_{T_\inda, E_{\inda-1}}^\ell-d_{E_{\inda}}^\ell\varphi_{T_\inda, F_\inda}^\ell)
 +\sum_{r=1}^{N_\kp}d_{T_\inda}^r\varphi_{T_\inda}^r\right)%
    \end{align}
    for arbitrary $d_0, d_{E_\inda}^\ell, d_{T_\inda}^r\in\R$ with $\ell=1,...,\kp, r=1,...,N_\kp, \inda=1,...,N$ (and $d_{E_0}^\ell=d_{E_{N}}^\ell$ for $z\in\mathcal V(\Omega)$) and the enumeration of $\T(z)$ as in Figure \ref{fig:NodePatch}.
    Furthermore, $\RzhD=\tszD + \szD$ holds for the solution
 $\szD\in \Vz$
to the $1+\kp|\F(z)|+\kp(\kp-1)/2N$-dimensional minimization problem \eqref{eqn:linear_min}.
\end{lemma}

\subsection{Proofs}\label{sub:app:proofs}
The remaining parts of this appendix are devoted to the verification of Lemmas \ref{lem:appendix}--\ref{lem:appendixB}.
\begin{proof}[of Lemma \ref{lem:appendix}]
Enumerate $\T(z)$ as in \ref{sub:app:alg} and
recall the definition of the jump $[\GG]_E = \GG|_{T_+} - \GG|_{T_-}$ on the interior edge $E=T_+\cap T_-$ shared by $T_+, T_-\in\T$, and $[\GG]_E = \GG|_{T_+}$ for the unique triangle $T_+\in\T$ with $E\subset T_+$ for the boundary edge $E\in\F(\partial \Omega)$. 
First, observe that $\sz\in\RTpwk(\mathcal{T}(z))$ lies in $\RTkpwz(\T(z))$ if and only if the coefficients $\cFj=0$ in the representation \eqref{eqn:basis_RTkpwz} are zero for $0\leq j\leq \kp$ at the edge $E\in\F(T)\setminus\F(z)$ in $T\in\T$ opposing $z$.
By definition, $\sz\in\RTkpwz(\T(z))$ belongs to $\Sz$ if and only if
\begin{align}
    \Ldivlm(\sz) &= (f+\ddiv \GG, \varphi_zx_1^\ell x_2^m)_{L^2(T)}&&\text{for all }0\leq \ell+m\leq \kp,T\in\T(z),\label{eqn:var_div}\\
    (\lambda_{T_+, E}^{j} + \lambda_{T_-, E}^{j})(\sz) &= ([\GG]_E\cdot \nu_E, \varphi_zs^j)_{L^2(E)}&&\text{for all }0\leq j\leq \kp, E\in \F(\omega(z)).\label{eqn:var_jump_cond}
\end{align}
This translates into equivalent conditions on the coefficients of $\sz$ in the representation \eqref{eqn:basis_RTkpwz}, namely, for all $\inda=1, ..., N$,
\begin{align}
    c_{T_\inda, E_{\inda}}^{0} &= (f+\ddiv \GG, \varphi_z)_{L^2(T_\inda)} - c_{T_\inda, E_{\inda-1}}^{0},&&\label{eqn:cond0}\\
    \cdivlma &= (f+\ddiv \GG, \varphi_zx_1^\ell x_2^m)_{L^2(T_\inda)}&&\text{for all }1\leq \ell+m\leq \kp,\label{eqn:c_div}\\
    c_{T_{\inda}, E_{\inda-1}}^j &= d_{E_{\inda-1}}^j&&\text{for all }0\leq j\leq \kp,\label{eqn:jump1}\\
    c_{T_\inda, E_\inda}^j &= ([\GG]_E\cdot \nu_E, \varphi_zs^j)_{L^2(E_\inda)} - d_{E_\inda}^j&&\text{for all }0\leq j\leq \kp,\label{eqn:jump2}
\end{align}
where $d_{E_\inda}^\ell\in\R$.
Since $\sz\in\RTkpwz(\T(z))$ vanishes at the other edges $E\in\F\setminus \F(z)$,  $\lambda_{T_\inda, \ddiv}^{0,0}(\sz) = c_{T_\inda, E_{\inda-1}}^{0} + c_{T_\inda, E_{\inda}}^{0}$ and \eqref{eqn:cond0}--\eqref{eqn:c_div} are equivalent to \eqref{eqn:var_div}.
The identification $d_{E_0}^\ell=d_{E_{N}}^\ell$ for an interior vertex $z\in\mathcal V(\Omega)$ with $E_0=E_N\in\F(z)$ shows that \eqref{eqn:jump1}--\eqref{eqn:jump2} are equivalent to \eqref{eqn:var_jump_cond}.
This identification is well defined.
Note that, whereas there is no condition on $d_{E_\inda}^\ell$ for $1\leq\ell\leq \kp$, the combination of \eqref{eqn:cond0} and \eqref{eqn:jump1} with \eqref{eqn:jump2} shows the implicit extra condition
$$d_{E_\inda}^0 = d_{E_0}^0 + \sum_{\alpha=1}^\inda \left(([\GG]_E\cdot \nu_E, \varphi_z)_{L^2(E_\alpha)} - (f+\ddiv \GG, \varphi_z)_{L^2(T_\alpha)}\right)\quad\text{for }\inda=0,...,N.$$
For an interior vertex $z\in\mathcal V(\Omega)$, an integration by parts and \eqref{eqn:crutial_identity} show that the sum on the right-hand side above vanishes for $a=N$, whence $d_{E_N}^0 = d_{E_0}^0$ is indeed well defined.
Furthermore, there is no condition on the coefficients $c_{T_\inda}^{r}$ for all $T_\inda\in\T(z)$ and $r=1,...,N_\kp$ and $c_{T_\inda}^{r}=d_{T_\inda}^r$ is a further degree of freedom.

Algorithm \ref{alg:eq} finds coefficients that satisfy \eqref{eqn:cond0}--\eqref{eqn:jump2} in a loop over $a=1,...,N$ and therefore defines $\tszD\in\Sz$ by \eqref{eqn:basis_RTkpwz}.\qed
\end{proof}
\begin{proof}[of Lemma \ref{lem:appendixB}]
This follows immediately after revisiting the proof of Lemma \ref{lem:appendix} for an arbitrary function $\sz\in \RTkpwz(\T(z))$.
Since $\sz\in \Sz$ is equivalent to \eqref{eqn:cond0}--\eqref{eqn:jump2} for the representation \eqref{eqn:basis_RTkpwz} of $\sz$ in the given basis, all functions $\sz\in \Sz -\tszD$ are of the form \eqref{eqn:sigma}
for arbitrary $d_0, d_{E_\inda}^\ell, d_{T_\inda}^r\in\R$ with $\ell=1,...,\kp$, $r=1,...,N_\kp$, and $\inda=1,...,N$ (and $d_{E_0}^\ell=d_{E_{N}}^\ell$ for $z\in\mathcal V(\Omega)$).
Hence, the dimension of the linear space $\Vz=\Sz -\tszD$ is $1+\kp|\F(z)|+\kp(\kp-1)/2N$.
The claim follows from Proposition \ref{prop:RzhD_def} by observing 
\begin{align*}
    \RzhD
\coloneqq
& 
\argmin\limits_{%
       \sz\in \Sz
     }
 \| \sz \|\niLz = \tszD + \argmin\limits_{%
       \sz\in \Sz - \tszD
     }
 \| \tszD + \sz \|\niLz.\quad\qed
\end{align*}
\end{proof}


%
%

\end{document}